\documentclass[12pt]{amsart}
\usepackage{amscd,amssymb,epsfig}
\usepackage[mathscr]{eucal}
\usepackage[english]{babel}
\usepackage{psfrag}

\input xy
\usepackage[all]{xy}

\begin{document}


\newtheoremstyle{mytheorem}
  {}
  {}
  {\slshape}
  {}
  {\scshape}
  {.}
  { }
  {}

\newtheoremstyle{mydefinition}
  {}
  {}
  {\upshape}
  {}
  {\scshape}
  {.}
  { }
  {}

\theoremstyle{mytheorem}
\newtheorem{lemma}{Lemma}[section]
\newtheorem{prop}[lemma]{Proposition}
\newtheorem{prop_intro}{Proposition}
\newtheorem{cor}[lemma]{Corollary}
\newtheorem{cor_intro}[prop_intro]{Corollary}
\newtheorem{thm}[lemma]{Theorem}
\newtheorem{thm_intro}[prop_intro]{Theorem}
\newtheorem*{thm*}{Theorem}
\theoremstyle{mydefinition}
\newtheorem{rem}[lemma]{Remark}
\newtheorem*{claim*}{Claim}
\newtheorem{rem_intro}[prop_intro]{Remark}
\newtheorem{rems_intro}[prop_intro]{Remarks}
\newtheorem*{notation*}{Notation}
\newtheorem*{warning*}{Warning}
\newtheorem{rems}[lemma]{Remarks}
\newtheorem{defi}[lemma]{Definition}
\newtheorem*{defi*}{Definition}
\newtheorem{defi_intro}[prop_intro]{Definition}
\newtheorem{defis}[lemma]{Definitions}
\newtheorem{exo}[lemma]{Example}
\newtheorem{exo_intro}[prop_intro]{Example}
\newtheorem{exos_intro}[prop_intro]{Examples}
\newtheorem*{exo*}{Example}
\newtheorem*{que*}{Question}

\numberwithin{equation}{section}

\newcommand{\bqn}{\begin{eqnarray*}}
\newcommand{\eqn}{\end{eqnarray*}}
\newcommand{\bq}{\begin{eqnarray}}
\newcommand{\eq}{\end{eqnarray}}
\newcommand{\ba}{\begin{aligned}}
\newcommand{\ea}{\end{aligned}}
\newcommand{\be}{\begin{enumerate}}
\newcommand{\ee}{\end{enumerate}}

\newcommand{\bibURL}[1]{{\unskip\nobreak\hfil\penalty50{\tt#1}}}

\def\ti{-\allowhyphens}

\newcommand{\thismonth}{\ifcase\month 
  \or January\or February\or March\or April\or May\or June%
  \or July\or August\or September\or October\or November%
  \or December\fi}
\newcommand{\thismonthyear}{{\thismonth} {\number\year}}
\newcommand{\thisdaymonthyear}{{\number\day} {\thismonth} {\number\year}}

\newcommand{\CC}{{\mathbb C}}
\newcommand{\DD}{{\mathbb D}}
\newcommand{\FF}{{\mathbb F}}
\newcommand{\HH}{{\mathbb H}}
\newcommand{\GG}{{\mathbb G}}
\newcommand{\KK}{{\mathbb K}}
\newcommand{\LL}{{\mathbb L}}
\newcommand{\NN}{{\mathbb N}}
\newcommand{\PP}{{\mathbb P}}
\newcommand{\QQ}{{\mathbb Q}}
\newcommand{\RR}{{\mathbb R}}
\renewcommand{\SS}{{\mathbb S}}
\newcommand{\TT}{{\mathbb T}}
\newcommand{\ZZ}{{\mathbb Z}}

\newcommand{\Bb}{{\mathcal B}}
\newcommand{\Cc}{{\mathcal C}}
\newcommand{\Dd}{{\mathcal D}}
\newcommand{\Ff}{{\mathcal F}}
\newcommand{\Ee}{{\mathcal E}}
\newcommand{\Gg}{{\mathcal G}}
\newcommand{\Hh}{{\mathcal H}}
\newcommand{\Kk}{{\mathcal K}}
\newcommand{\Ll}{{\mathcal L}}
\newcommand{\Mm}{{\mathcal M}}
\newcommand{\Nn}{{\mathcal N}}
\newcommand{\Oo}{{\mathcal O}}
\newcommand{\Pp}{{\mathcal P}}
\newcommand{\Qq}{{\mathcal Q}}
\newcommand{\Rr}{{\mathcal R}}
\newcommand{\Ss}{{\mathcal S}}
\newcommand{\Tt}{{\mathcal T}}
\newcommand{\Vv}{{\mathcal V}}
\newcommand{\Xx}{{\mathcal X}}
\newcommand{\Yy}{{\mathcal Y}}
\newcommand{\Zz}{{\mathcal Z}}

\newcommand{\fraka}{{\mathfrak a}}
\newcommand{\frakg}{{\mathfrak g}}
\newcommand{\frakk}{{\mathfrak k}}
\newcommand{\frakl}{{\mathfrak l}}
\newcommand{\frakp}{{\mathfrak p}}
\newcommand{\fraks}{{\mathfrak s}}
\newcommand{\fraku}{{\mathfrak u}}

\newcommand{\dD}{{\mathbf D}}
\newcommand{\gG}{{\mathbf G}}
\newcommand{\hH}{{\mathbf H}}
\newcommand{\pP}{{\mathbf P}}
\newcommand{\lL}{{\mathbf L}}
\newcommand{\qQ}{{\mathbf Q}}

\newcommand{\<}{\langle}
\renewcommand{\>}{\rangle}

\def\binfty{\mathcal B^\infty_{\mathrm {alt}}}
\def\bu{\bullet}
\def\cb{{\rm C}_{\rm b}}
\def\ehbc{{\rm EH}_{\rm cb}}
\def\Ef{\mathcal E_\varphi}
\def\essim{\operatorname{EssIm}}
\newcommand{\essimfi}{\operatorname{EssIm(\varphi)}}
\newcommand{\esssup}{\operatorname{ess\,sup}}
\def\h{{\rm H}}
\def\hb{{\rm H}_{\rm b}}
\def\hc{{\rm H}_{\rm c}}
\def\hcb{{\rm H}_{\rm cb}}
\def\hom{\operatorname{Hom}}
\def\homeo#1{{\sl H\!omeo}^+\!\left(#1\right)}
\def\thomeo#1{\widetilde{{\sl \!H}\!omeo}^+\!\left(#1\right)}
\def\id{{\it I\! d}}
\def\ind{\mathrm{ind}}
\def\la{\mathrm{L}^\infty_{\mathrm{alt}}}
\def\linfty{\mathrm{L}^\infty}
\def\linftyw{\mathrm{L}^\infty_{\rm w*}}
\def\lp{\mathrm{L}^p}
\def\ltwo{\mathrm{L}^2}
\def\rep{\operatorname{Rep}}
\newcommand{\supp}{\operatorname{supp}}

\def\one{\mathbf{1\kern-1.6mm 1}}
\def\property{\textbf{\rm\textbf A}}
\def\oddex#1#2{\left\{#1\right\}_{o}^{#2}}
\def\comp#1{{\rm C}^{(#1)}}
\def\lra{\longrightarrow}

\def\adg{\operatorname{ad}_\frakg}
\def\bg{B_\frakg}
\def\cs{\check S}
\def\deta{{\operatorname{det}_A}}
\def\gl{\mathrm{GL}}
\def\gg{\Gamma_g}
\def\ghtp{generalized Hermitian triple product }
\def\gmodp{\gG(\RR)/\pP(\RR)}
\def\gmodq{\gG(\RR)/\qQ(\RR)}
\def\gr{\mathrm{Gr}_p(W)}
\def\gra{\mathrm{Gr}^A_p(V)}
\def\Gr{\mathrm{Gr}}
\def\htp{Hermitian triple product }
\newcommand{\Iso}{\operatorname{Iso}}
\def\Isom{\operatorname{Isom}}
\def\isp{\operatorname{Is}_{\<\cdot,\cdot\>}}
\def\isptwo{\operatorname{Is}_{\<\cdot,\cdot\>}^{(2)}}
\def\ispth{\operatorname{Is}_{\<\cdot,\cdot\>}^{(3)}}
\def\isf{\operatorname{Is}_F}
\def\isfi{\operatorname{Is}_{F_i}}
\def\isft{\operatorname{Is}_F^{(3)}}
\def\isfit{\operatorname{Is}_{F_i}^{(3)}}
\def\isftwo{\operatorname{Is}_F^{(2)}}
\def\lin{\operatorname{Lin}(L_+,L_-)}
\def\ll{{\Ll_1,\Ll_2}}
\def\kahler{K\"ahler }
\def\kg{\kappa_G}
\def\kgb{\kappa_G^{\rm b}}
\def\kib{\kappa_i^{\rm b}}
\def\kibt{\tilde\kappa_i^{\rm b}}
\def\kx{\kappa_\Xx}
\def\kxb{\kappa_\Xx^\mathrm{b}}
\def\ox{\omega_\Xx}
\def\po{\mathrm{PO}}
\def\pr{\operatorname{pr}}
\def\psl{\mathrm{PSL}}
\def\pu{\mathrm{PU}}
\def\pupq{\mathrm{PU}(p,q)}
\def\puvi{\mathrm{PU}\big(V,\<\cdot,\cdot\>_i\big)}
\def\rkx{\operatorname{rk}_\Xx}
\def\rky{\operatorname{rk}_\Yy}
\def\sg{\Sigma_g}
\def\sltwo{\mathrm{SL}(2,\RR)}
\def\slv{\mathrm{SL}(V)}
\def\sp{\mathrm{Sp}}
\def\stab{\operatorname{Stab}}
\def\supq{\mathrm{SU}(p,q)}
\def\su{\mathrm{SU}}
\def\suq{\mathrm{SU}(q,1)}
\def\suw{\mathrm{SU}(W)}
\def\suvi{{\rm SU}\big(V,\<\cdot,\cdot\>_i\big)}
\def\tr{\mathrm{T}_\rho}
\def\u{\mathrm{U}}
\def\vol{\operatorname{vol}}
\def\bsl{\backslash}


\title[Maximal Representations of Surface Groups]{Maximal Representations\\ of Surface Groups:\\Symplectic Anosov Structures}
\author[M.~Burger]{Marc Burger}
\email{burger@math.ethz.ch}
\address{FIM, ETH Zentrum, R\"amistrasse 101, CH-8092 Z\"urich, Switzerland}
\author[A.~Iozzi]{Alessandra Iozzi}
\email{Alessandra.Iozzi@unibas.ch}
\address{Institut f\"ur Mathematik, Universit\"at Basel, Rheinsprung 21,
CH-4051 Basel, Switzerland}
\address{D\'epartment de Math\'ematiques, Universit\'e de Strasbourg, 7, rue Ren\'e Descartes, F-67084 Strasbourg Cedex, France}
\thanks{A.I. and A.W. were partially supported by FNS grant PP002-102765}
\author[F.~Labourie]{Fran\c cois Labourie}
\email{labourie@math.u-psud.fr}
\address{Topologie et Dynamique, B\^atiment 425, Universit\'e de Paris Sud, F-91405 Orsay Cedex, France}
\author[A.~Wienhard]{Anna Wienhard}
\email{Anna.Wienhard@unibas.ch}
\address{Institut f\"ur Mathematik, Universit\"at Basel, Rheinsprung 21,
CH-4051 Basel, Switzerland}

\date{\today}

\maketitle

\centerline{\it Dedicated to the memory of Armand Borel,}
\centerline{\it with affection and admiration}
\vskip2cm

\begin{abstract}
Let $G$ be a connected semisimple Lie group such that the associated symmetric
space $\Xx$ is Hermitian and let $\gg$ be the fundamental group of 
a compact orientable surface of genus $g\geq2$.
We survey the study of maximal representations of $\gg$ into $G$,
that is the subset of $\hom(\gg,G)$  
characterized by the maximality of the Toledo invariant 
(\cite{Burger_Iozzi_Wienhard_ann} and \cite{Burger_Iozzi_Wienhard_tol}).  
Then we concentrate on the particular case $G=\sp(2n,\RR)$,
and we show that if $\rho$ is any maximal representation 
then the image $\rho(\gg)$ is a discrete, faithful realizations of $\gg$ 
as a Kleinian group of complex motions in $\Xx$ with an associated Anosov system,
and whose limit set in an appropriate compactification of $\Xx$ 
is a rectifiable circle.
\end{abstract}

\vskip2cm

\setcounter{tocdepth}{1}
\tableofcontents

\section{Introduction}\label{sec:intro}
Let $\gg$ be the fundamental group of a compact orientable topological surface 
$\sg$ of genus $g\geq2$.  
For a general real algebraic group $G$ the representation variety $\hom(\gg,G)$
is a natural geometric object which reflects properties both of the discrete
group $\gg$ and of the algebraic group $G$ and enjoys an extremely rich
structure.  For example, $\hom(\gg,G)$ is not only a topological space,
but also a real algebraic variety, which in addition 
parametrizes flat principal $G$-bundles over $\sg$;
furthermore, it admits an action of the group of automorphisms of $\gg$
by precomposition which commutes with the action by postcomposition
with (inner) automorphisms of $G$. 
It is natural to consider homomorphisms up to conjugation,
thus we introduce the topological quotient
\bqn
\rep(\gg,G):=\hom(\gg,G)/G\,;
\eqn

although this quotient is not necessarily a Hausdorff space,
it contains a large part which is Hausdorff,
namely the space $\rep_{red}(\gg,G)$ of homomorphisms
with reductive image modded out by $G$-conjugation.
The general theme of this note is the study of certain
connected components of $\hom(\gg,G)$ or $\rep(\gg,G)$
analogous to Teichm\"uller space,
and their relation to geometric and dynamical structures on $\sg$.  

Recall that if $G=\pu(1,1)$, the space $\rep(\gg,G)$ 
has $4g-3$ connected components (\cite{Goldman_thesis}, \cite{Goldman_88}), 
two of which
are homeomorphic to $\RR^{6g-6}$ and correspond to the two Teichm\"uller spaces 
$\Tt_g$ -- one for each orientation of $\sg$ -- 
that is to the space of marked complex, alternatively hyperbolic, 
structures on the topological surface $\sg$.

If on the other hand $G=\mathrm{SL}(3,\RR)$, 
Goldman and Choi proved \cite{Goldman_Choi}
that one of the three connected components of $\hom(\gg,G)$ \cite{Hitchin}
parametrizes convex projective structures on $\sg$,
that is diffeomorphisms of $\sg$ with $\Omega/\Gamma$, where 
$\Gamma<\mathrm{SL}(3,\RR)$ is a faithful discrete image of $\gg$
and $\Omega\subset\PP(\RR^3)$ is a convex invariant domain;
incidentally, this component coincides with the Hitchin component
that we define below.

If $G=\psl(2,\CC)$, there is an open subset of $\rep(\gg,G)$
consisting of all quasi-Fuchsian deformations of $\gg$,
which is diffeomorphic to the product $\Tt_g\times\Tt_g$ 
of two copies of Teichm\"uller space.

In each of these three cases, a representation $\rho$
belonging to such a ``special component'' in $\rep(\gg,G)$
is faithful and with discrete image, and $\rho(\gg)<G$
gives rise, as a Kleinian group, to many interesting dynamical
and geometric structures.

When $G$ is a simple split real Lie group --
such as for instance $G=\psl(n,\RR)$, $\mathrm{PSp}(2n,\RR)$,
$\mathrm{PO}(n,n)$ or $\mathrm{PO}(n,n+1)$ -- 
Hitchin singled out a component $\rep_H(\gg,G)$ of $\rep(\gg,G)$
which he proved, using Higgs bundle techniques, 
to be diffeomorphic to $\RR^{|\chi(\sg)|\dim G}$, \cite{Hitchin}
and which is now commonly known as Hitchin component.
For example, if $G=\psl(n,\RR)$ for $n\geq2$, 
this component is the one containing the homomorphisms
of $\gg$ into $\mathrm{SL}(n,\RR)$ obtained by composing
a hyperbolization with the $n$-dimensional irreducible
representation of $\sltwo$.   As Hitchin however points out
\cite{Hitchin},
the analytic point of view does not shed any light on the
geometric significance of the representations in this component.

Recently the concept of Anosov representation,
which links the surface $\sg$ to flag manifolds associated to $G$
was introduced in \cite{Labourie_anosov} 
and used to show that, if $G=\psl(n,\RR)$, representations in the Hitchin component
are discrete and faithful, and that they provide quasiisometric
embeddings of $\gg$ into $G$, \cite{Labourie_anosov}, \cite{Labourie_crossratio}.

In parallel, Goncharov and Fock developed for surfaces
with nonempty boundary and for split real 
Lie groups a tropical-algebro-geometric viewpoint of 
$\rep(\gg,G)$, singling out positive real points in $\rep(\gg,G)$ 
which correspond to discrete and faithful representations,
\cite{Fock_Goncharov}, \cite{Fock_Goncharov_04}.

There is another natural extension of the case $G=\pu(1,1)$ in a different
direction, that is to connected semisimple Lie groups $G$
such that the associated symmetric space $\Xx$ admits
a $G$-invariant complex structure, just like in the case of
the Poincar\'e disk.  This includes notably groups
like $\su(p,q)$, $\sp(2n,\RR)$, $\mathrm{SO}^\ast(2n)$,
$\mathrm{SO}(2,n)$.  Symmetric spaces with this property
are called {\it Hermitian}.

In the same framework, the topology and the number of connected
components of the space of reductive representations into $\su(p,q)$ 
and $\sp(2n,\RR)$ have been studied in a series of papers 
by Bradlow, Garc\'ia-Prada, Gothen, Mundet i Riera and Xia
(\cite{Bradlow_GarciaPrada_Gothen_big}, 
\cite{Bradlow_GarciaPrada_Gothen}, \cite{Bradlow_GarciaPrada_Gothen_cc},
\cite{Gothen}, \cite{GarciaPrada_Mundet}, \cite{Xia}), 
extending the analytic approach introduced by Hitchin.

The additional feature for symmetric spaces which are
Hermitian is the presence of a \kahler form $\ox$ on $\Xx$ 
which allows to associate to every representation 
$\rho:\gg\to G$ a characteristic number,
called the {\it Toledo invariant} $\tr$ (see \S~\ref{sec:inv}), 
which is constant on connected components of $\hom(\gg,G)$ and 
which satisfies a Milnor--Wood type inequality 
\bq\label{eq:mw}
|\tr|\leq|\chi(\sg)|\rkx\,,
\eq
where $\rkx$ is the real rank of $X$.  A representation 
is {\it maximal} if equality holds in (\ref{eq:mw}),
and the set $\hom_{max}(\gg,G)$ of such representations is
then a union of components of $\hom(\gg,G)$.
If $G=\pu(1,1)$, Goldman proved in \cite{Goldman_thesis} 
that maximal representations are exactly those lying in the 
two Teichm\"uller components.  

In the first part of this article we illustrate, mostly without proofs, 
results concerning the geometric significance of maximal representations.
To fix the notation, let $G:=\gG(\RR)^\circ$, 
where $\gG$ is a semisimple real algebraic group and 
assume that the symmetric space $\Xx$ associated to $G$ is Hermitian.
In complete analogy with Goldman's theorem,
any maximal representation $\rho:\gg\to G$ is injective 
with discrete image (Theorem~\ref{thm:main}).
This fact depends on a careful study of the Zariski closure $\lL$ of
the image $\rho(\gg)$ and the fact that there is an essential 
restriction on $L:=\lL(\RR)^\circ$, namely that it is reductive
and it preserves a subHermitian symmetric space of $\Xx$ 
which is of {\it tube type} and maximal with respect to this property.
On the constructive side, the study of maximal representations into $G$ 
does not reduce to the study of classical Teichm\"uller space; in fact, 
if $\Xx$ is of tube type, any representation which is the 
composition of a hyperbolization $\gg\to\su(1,1)$ with the homomorphism 
$\su(1,1)\to G$ associated to the realization of the Poincar\'e disk
diagonally in a maximal polydisk in $\Xx$ can be deformed into
a representation with Zariski dense image in $\gG$
(Theorem~\ref{thm:ex}); such a representation is by construction maximal.
For examples of discrete representations into $\su(1,n)$
with prescribed Toledo invariant see \cite{Goldman_Kapovich_Leeb}.

These results are proven in greater generality in \cite{Burger_Iozzi_Wienhard_tol}, where for the representation of the fundamental group of a
surface with boundary, we define a Toledo invariant 
whose definition and properties however require 
some vigorous use of bounded cohomology.
In the context of this paper, continuous bounded cohomology 
will appear as a tool in the proofs; 
in particular it allows to define the notion of {\it tight homomorphism},
more general and flexible than that of maximal representation, 
and which is an essential tool to study the geometric properties 
of the inclusion $\Xx_L\hookrightarrow\Xx$, where $\Xx_L$ 
is the subsymmetric space associated to $L$ (see above).
A systematic study of tight homomorphisms and the companion notion
of {\it tight embedding} is the subject matter of 
\cite{Burger_Iozzi_Wienhard_tight} and a few highlights of this theory are
presented in \S~\ref{sec:tight}.

While the first part of the paper is expository, in the second part
we give an elementary treatment of a certain  number of results 
on maximal representations into the symplectic group $\sp(V)$
of a real symplectic vector space $V$.
The results are stated in \S~\ref{sec:anosov} and their proofs
in \S~\ref{sec:anproofs} are independent of the rest of the paper
(see \S~\ref{sec:anproofs}).

Observe at this point that $\sp(V)$ is 
at the same time real split, and hence falls into the context 
of the Hitchin component, and is the group of automorphisms of
the Siegel upper half space, a fundamental class of Hermitian symmetric spaces.
We have the inclusion
\bqn
       \rep_H\big(\gg,\sp(V)\big)
\subset\rep_{max}\big(\gg,\sp(V)\big)\,,
\eqn
but while the representations in the Hitchin component are all 
irreducible \cite{Labourie_anosov}, 
there are (at least when $\dim V\geq4$) components of maximal 
representations which contain reducible representations, 
so that the above inclusion is strict.

For a representation $\rho:\gg\to\sp(V)$ and a fixed hyperbolization
$\Sigma$ of $\sg$, we associate the flat symplectic bundle $E^\rho$
over the unit tangent bundle $T^1\Sigma$ of $\Sigma$ with fiber $V$.
The geodesic flow lifts canonically to a flow $g_t^\rho$ on $E^\rho$;
we adapt some of the ideas in \cite{Labourie_anosov} to our
situation and, combining them with the results in \S~\ref{sec:tight} and 
\S~\ref{sec:proofs}, prove that if $\rho$ is maximal then $E^\rho$ 
is the sum of two continuous Lagrangian subbundles $E^\rho_+\oplus E^\rho_-$
on which $g_t^\rho$ acts contracting and expanding respectively.
Moreover, this bundle will also come with a field of complex
structures in each fiber, exchanging $E_\pm^\rho$ and 
positive for the symplectic structure (see \S~\ref{sec:tight}).  
As a consequence, one deduces that any maximal representation 
$\rho:\gg\to\sp(V)$ is a quasiisometric embedding,
where $\sp(V)$ is equipped with a standard invariant metric.
This implies that the action of the mapping class group
$\operatorname{Out}(\gg)$ on $\rep_{max}\big(\gg,\sp(V)\big)$
is properly discontinuous.

\vskip.7cm
{\it Acknowledgments:} The authors thank 
Domingo Toledo for bringing to their attention \cite{Kneser} and \cite{Toledo_email}
and Domingo Toledo and Nicolas Monod for useful comments 
on a preliminary version of the paper.


\vskip1cm

\section{Hermitian Symmetric Spaces and Examples}\label{sec:hss}
Let $\Xx$ be a symmetric space and let $G:=\Isom(\Xx)^\circ$ be 
the connected component of its group of isometries.
Recall that $\Xx$ is {\it Hermitian} if it admits
a $G$-invariant complex structure.  An equivalent definition
is that $\Xx$ is a Hermitian manifold such that every point $x\in\Xx$ 
is the isolated fixed point of an isometric involution $s_x$.
In this paper we shall consider only symmetric spaces
of noncompact type.

Let $J$ be the $G$-invariant complex structure
and let $g:T\Xx\times_p T\Xx\to\RR$ be the Riemannian metric, where 
$T\Xx\times_p T\Xx$ denotes the fibered product over the projection $p:T\Xx\to \Xx$. 
Then 
\bqn
\omega_\Xx(X,Y):=g(JX,Y)
\eqn
defines a $G$-invariant differential two-form on $\Xx$
which is nondegenerate.

\begin{lemma}\label{lem:closed} 
Let $\Xx$ be a symmetric space and $G=\Isom(\Xx)^\circ$.
Then any $G$-invariant differential form on $\Xx$ is closed.
\end{lemma}
\begin{proof}
Let $\alpha$ be a $G$-invariant differential $k$-form on $\Xx$ 
and let $s\in\Isom(\Xx)$ be the geodesic symmetry at a basepoint $0\in\Xx$.  
Since $G$ is normal in $\Isom(\Xx)$,
then $sgs\in G$ and hence $s\alpha$ is also $G$-invariant, since
\bqn
g(s\alpha)=s^2g(s\alpha)=s(sgs)\alpha=s\alpha\,.
\eqn 
Moreover, since $s|_{T_0\Xx}=-\id$ we have that 
$(s\alpha)_0=(-1)^k\alpha_0$, and since $\alpha$ and $s\alpha$ 
are both $G$-invariant, the equality $(s\alpha)_x=(-1)^k\alpha_x$ 
holds for every $x\in\Xx$. Since $d\alpha$ is also $G$-invariant, from
\bqn
(-1)^kd\alpha=d(s\alpha)=s(d\alpha)=(-1)^{k+1}d\alpha
\eqn
we deduce that $d\alpha=0$.
\end{proof}

The immediate consequence of the above lemma is that a Hermitian
symmetric space $\Xx$ is a \kahler manifold with \kahler form $\omega_\Xx$.
Furthermore, using the existence of a \kahler form
on $\Xx$, one can prove that for an irreducible symmetric space $\Xx$
being Hermitian is equivalent to the center of a maximal compact
subgroup of $\Isom(\Xx)^\circ$ having positive dimension 
(and in fact being one-dimensional).

\medskip

A fundamental result which makes the study of Hermitian
symmetric spaces accessible to techniques from function theory 
{\it \`a la Bergmann} is the following theorem of Harish-Chandra
which for classical domains is due to E.~Cartan, \cite{Cartan_35}. 

\begin{thm}[Harish-Chandra, \cite{Harish-Chandra_VI}] 
Any Hermitian symmetric space of noncompact type
is biholomorphic to a bounded domain in a complex vector space.
\end{thm}

The bounded realization $\Dd\subset\CC^N$ of a Hermitian symmetric space $\Xx$
has a natural compactification, namely the topological closure
$\overline\Dd$ in $\CC^N$, on which 
$G:=\Isom(\Xx)^\circ$ acts by restriction of birational isomorphism of $\CC^N$. 
The {\it Shilov boundary} $\cs$ is a subset of the topological boundary 
$\partial\Dd$ of the bounded domain; 
it can be defined in function theoretical terms, and it is also the unique closed
$G$-orbit in $\overline\Dd$.  It is a homogeneous space 
of the form $G/Q$, where $Q$ is a (specific) maximal parabolic subgroup of $G$,
and plays a prominent role in our study, for example as target of
appropriate boundary maps.  Notice that only if $\Xx$ is of real 
rank one, the Shilov boundary coincides with the whole boundary $\partial\Dd$.

Recall that the rank $\rkx$ of a symmetric space $\Xx$ 
is the maximal dimension of a flat subspace,
that is an isometric copy of Euclidean space.

Expositions of different aspects of the geometry of Hermitian symmetric spaces
are \cite{Koranyi_infarautetal}, \cite{Deligne}, \cite{Satake_book},
\cite{Piateskii}, and \cite{Wolf}.

\subsection{Examples of Hermitian Symmetric Spaces}
We give here examples of two families of Hermitian symmetric spaces
which are of fundamental nature and with which we shall
illustrate our results.

\subsubsection{$\suw$}\label{subsubsec:suw} 
Let $W$ be a complex vector space of dimension $n$,
and $h(\,\cdot\,,\,\cdot\,)$ a nondegenerate Hermitian form of signature
$(p,q)$, $p\leq q$, so that $p$ is the maximal dimension of a subspace $L\subset W$
on which the restriction $h(\,\cdot\,,\,\cdot\,)|_L$ is positive
definite.  A model for the symmetric space associated to 
\bqn
\suw:=\big\{g\in\mathrm{SL}(W):\, \,h(gx,gy)=h(x,y), \forall x,y\in W\big\}
\eqn is
\bqn
\Xx_{\suw}:=\big\{L\in\gr: h(\,\cdot\,,\,\cdot\,)|_L\hbox{ is positive definite}\big\}
\eqn
which, as an open subset of the Grassmannian $\gr$ of $p$-dimensional subspaces of $W$, 
is a complex manifold on which $G$ acts by automorphisms.  
To realize $\Xx_{\suw}$ as a bounded domain, 
fix $W_+\in\Xx_{\suw}$, and let $W_-:=W_+^\bot$ be its orthogonal complement
with respect to the form $h$.  Since the orthogonal
projection ${\rm pr}_\epsilon:W\to W_\epsilon$, $\epsilon\in\{+,-\}$, 
is an isomorphism for $\epsilon=+$ when restricted
to any $L\in \Xx_{\suw}$, we can define 
\bq\label{eq:E}
E:\Xx_{\suw}\to&\operatorname{Lin}(W_+,W_-)
\eq
by
\bqn
E(L):={\rm pr}_-\circ({\rm pr}_+|_L)^{-1}\,.
\eqn
It is easy to see that this defines a biholomorphic map from $\Xx_{\suw}$ 
to the bounded domain
\begin{equation}\label{eq:dpq}
\Dd_{\suw}:=\big\{A\in\operatorname{Lin}(W_+,W_-):
\,\id-A^\ast A\hbox{ is positive definite}\big\}\,
\end{equation}
where the adjoint is taken with respect to the structures of 
the unitary spaces $\big(W_+,h(\,\cdot\,,\,\cdot\,)|_{W_+}\big)$ and 
$\big(W_-,-h(\,\cdot\,,\,\cdot\,)|_{W_-}\big)$.  Moreover the closure
$\overline\Xx_{\suw}$ of $\Xx_{\suw}$ in $\gr$ is mapped by $E$ to 
\bqn
\overline{\Dd_{\suw}}=\big\{A\in\operatorname{Lin}(W_+,W_-):\id-A^\ast A\hbox{ is positive semidefinite}\big\}\,.
\eqn
To determine the preimage of the Shilov boundary in the hyperboloid model
$\Xx_{\suw}$, observe that there are precisely $(p+1)$ orbits of $\suw$ 
in $\overline\Xx_{\suw}$, 
only one of which is closed, namely the Grassmannian of maximal isotropic subspaces
\bqn
\operatorname{Is}_p(W):=\big\{L\in\gr:\,h(\,\cdot\,,\,\cdot\,)|_L=0\big\}\,,
\eqn
which is hence sent via $E$ to the Shilov boundary 
\bqn
\cs_{\suw}=\big\{A\in\operatorname{Lin}(W_+,W_-):\,\id-A^\ast A=0\big\}
\subset\overline{\Dd_{\suw}}
\eqn
of the bounded domain $\Dd_{\suw}$.  
The real rank of $\Xx_{\suw}$ is $p$.  

Identifying $W$ with $\CC^{p+q}$ in such a way that
$h$ is the standard Hermitian form of signature $(p,q)$,
we denote $\supq:=\suw$ and
\bqn
\Dd_{p,q}=\big\{Z\in M_{q,p}(\CC):\id-{}^t\overline{Z}Z\hbox{ is positive definite}\big\}
\eqn
the corresponding bounded domain with Shilov boundary
\bqn
\cs_{p,q}=\big\{Z\in M_{q,p}(\CC):\id-{}^t\overline{Z}Z=0\big\}\subset
\overline{\Dd_{p,q}}\,.
\eqn
In particular $\Dd_{1,1}$ is the Poincar\'e disk.

\subsubsection{The Symplectic Group}\label{subsubsec:sympl}  
Let $V$ be a real vector space equipped with a symplectic form
$\<\,\cdot\,,\,\cdot\,\>$, that is a nondegenerate antisymmetric bilinear form.
In particular $V$ must be even dimensional and we fix $\dim V=2n$.
The group 
\bqn
\sp(V):=\big\{g\in\gl(V):\, \,\<gx,gy\>=\<x,y\>, \forall x,y\in V\big\}
\eqn
is the real symplectic group.  The fact that on a complex vector space 
the imaginary part of a nondegenerate Hermitian form is a symplectic form
for the underlying real structure suggests to introduce the space
\bqn
\begin{aligned}
\Xx:=\big\{ J\in\gl(V):&\, J\hbox{ is a complex structure on $V$ and }\\
&h_J(x,y):=\<x,Jy\>+i\<x,y\>\hbox{ is a positive definite}\\
&\hbox{Hermitian form on }(V,J)\big\}\,,
\end{aligned}
\eqn
so that, if $J\in\Xx$, then
$\Re h_J(x,y)=\<x,Jy\>$ is a symmetric positive definite form.
It is easy to see that, among complex structures on $V$, 
this property characterizes the elements of $\Xx$.
Furthermore, for the transitive action by conjugation of
$\sp(V)$ on $\Xx$, the stabilizer of $J$ is a maximal compact subgroup
isomorphic to $\u(n)$ and hence $\Xx$ is the symmetric space $\Xx_{\sp(V)}$ 
associated to $\sp(V)$; 
in particular, since the center of $\u(n)$ has positive dimension,
$\Xx_{\sp(V)}$ is Hermitian symmetric and, as such,
there is a $\sp(V)$-invariant complex structure on $\Xx_{\sp(V)}$
which one can explicit as follows.

Let $V_\CC$ be the complexification of $V$ 
and let $\sigma:V_\CC\to V_\CC$ be the complex conjugation $\sigma(x+iy):=x-iy$
for $x,y\in V$.
Then there is a bijective correspondence between complex structures $J$ on $V$ and 
decompositions $V_\CC=L_+\oplus L_-$ into complex subspaces
satisfying $\sigma(L_\pm)=L_\mp$, given by the eigenspace
decomposition of $J\otimes_\CC\one$ into $\pm i$-eigenspaces.  
If now $\<\,\cdot\,,\,\cdot\,\>_\CC$ is the complexification of 
$\<\,\cdot\,,\,\cdot\,\>$, then
\bqn
h(x,y):=i\<x,\sigma(y)\>_\CC
\eqn
is a nondegenerate Hermitian form of signature $(n,n)$ on $V_\CC$;
if in particular $J\in\Xx_{\sp(V)}$, then
the restriction $h|_{L_+\times L_+}$ is positive definite, 
so that we obtain a map
\bq\label{eq:embedding}
\ba
\Xx_{\sp(V)}&\to\Xx_{\su(V_\CC)}\\
J\quad&\mapsto \quad L_+
\ea
\eq
which is equivariant with respect to the natural homomorphism 
\bqn
\begin{aligned}
\lambda:\sp(V)&\to\su(V_\CC)\\
g\quad&\mapsto g\otimes_\CC\one\,.
\end{aligned}
\eqn
Since $\Xx_{\su(V_\CC)}$ inherits a natural complex structure 
as an open subset of the Grassmannian 
$\Gr_n(V_\CC)$ of $n$-planes in $V_\CC$ and
\bqn
M:=\big\{L\in\Gr_n(V_\CC):\,\<\,\cdot\,,\,\cdot\,\>_\CC|_{L\times L}\equiv0\big\}
\eqn
is an algebraic subvariety of $\Gr_n(V_\CC)$, then
$\Xx_{\su(V_\CC)}\cap M$ acquires a natural complex structure
as an open subset of $M$. But  
\bqn
\<\,\cdot\,,\,\cdot\,\>_\CC|_{L_\pm\times L_\pm}\equiv0
\eqn
for all $J\in\Xx_{\sp(V)}$,
so that the map in (\ref{eq:embedding}) is 
a $\lambda$-equivariant bijection between $\Xx_{\sp(V)}$ and 
$\Xx_{\su(V_\CC)}\cap M$ by the use of which 
the complex structure on $\Xx_{\su(V_\CC)}\cap M$ defines 
the $\sp(V)$-invariant complex structure on $\Xx_{\sp(V)}$.

Let us denote by $\Ss_{\sp(V)}:=\Xx_{\su(V_\CC)}\cap M$ 
the {\it Siegel space} associated to $\sp(V)$
with its $\sp(V)$-action via the homomorphism $\lambda$.  

To determine the bounded domain realization of $\Xx_{\sp(V)}$, 
it is enough to observe that 
-- with the notations of \S~\ref{subsubsec:suw}, where $W=V_\CC$ -- 
the Siegel space $\Ss_{\sp(V)}$ is mapped
by the map $E$ defined in (\ref{eq:E}) to the subdomain
\bqn
\Dd_{\sp(V)}:=\big\{
A\in\Dd_{\su(V_\CC)}:\,\<\,\cdot\,,A\,\cdot\,\>_\CC|_{W_+\times W_+}
\hbox{ is symmetric}\big\}\,,
\eqn
and, accordingly, $\overline\Ss_{\sp(V)}$ is mapped to 
\bqn
\overline{\Dd_{\sp(V)}}:=\big\{
A\in\overline{\Dd_{\su(V_\CC)}}:\,\,\<\,\cdot\,,A\,\cdot\,\>_\CC|_{W_+\times W_+}
\hbox{ is symmetric}\big\}\,.
\eqn
One can verify again that the closure $\overline{\Ss_{\sp(V)}}$ in 
$\Gr_n(V_\CC)$ decomposes into $(n+1)$ orbits under the symplectic group,
only one of which is closed, namely
\bqn
 \mathrm{Is}_n(V_\CC)\cap\overline{\Ss_{\sp(V)}}
=\big\{L\in\Gr_n(V_\CC):\,
\<\,\cdot\,,\,\cdot\,\>_\CC|_{L\times L}=0\hbox{ and }\,h|_{L\times L}\equiv0 \big\}
\eqn
and hence is the preimage, under the map $E$ in (\ref{eq:E}), 
of the Shilov boundary in the bounded domain realization of $\Ss_{\sp(V)}$. 
To give an intrinsic description of the Shilov boundary 
observe that, since we have the alternative description
\bqn
\mathrm{Is}_n(V_\CC)\cap\overline{\Ss_{\sp(V)}}
=\big\{
L\in\Gr_n(V_\CC):\,\<\,\cdot\,,\,\cdot\,\>_\CC|_{L\times L}=0\hbox{ and }\sigma(L)=L \big\}\,,
\eqn
we conclude that the map
\bqn
\begin{aligned}
\Ll(V)&\to\mathrm{Is}_n(V_\CC)\cap\overline{\Ss_{\sp(V)}}\\
L\quad&\mapsto \quad L\otimes\CC\,,
\end{aligned}
\eqn
where $\Ll(V)$ is the space of Lagrangians in $V$,
is an $\sp(V)$-equivariant bijection.
Thus the space of Lagrangian subspaces is, as a $\sp(V)$-homogeneous 
space, isomorphic to the Shilov boundary of the bounded domain $\Dd_{\sp(V)}$.

If we identify $V$ with the direct sum of $n$ symplectic planes,
that is copies of $\RR^2$ with the standard symplectic form,
then accordingly we denote the symplectic group by $\sp(2n,\RR)$ and
the Siegel space by $\Ss_n$.


\vskip1cm

\section{The Toledo Invariant and Maximal Representations}\label{sec:inv}
Let $\Sigma_g$ be a compact oriented surface of genus $g\geq2$,
$G$ a connected semisimple Lie group with finite center and associated
symmetric space $\Xx$ which we assume to be Hermitian, 
and $\rho:\gg\to G$ a homomorphism.  
Then there is a smooth $\gg$-equivariant map $\tilde f:\widetilde\sg\to\Xx$, 
where $\widetilde\sg$ denotes the universal covering of $\sg$,
which can be obtained by lifting a smooth section of the flat bundle
$\gg\backslash\big(\widetilde{\sg}\times\Xx\big)\to\sg$ with contractible fiber $\Xx$.
The pullback $\tilde f^\ast\omega_\Xx$
is then a $\Gamma$-invariant differential form on $\widetilde\sg$, 
which hence descends to a form on $\sg$.
Since any two such sections are homotopic
and hence the map $\tilde f$ is unique up to $\rho$-equivariant
homotopy, the result of the integration 
over $\sg$ of any two forms obtained in this way
does not depend on the particular choice of a section: we can hence
define the {\it Toledo invariant of $\rho$}
\bqn
\tr:=\frac{1}{2\pi}\int_{\sg}\tilde f^\ast(\omega_\Xx)\,.
\eqn

Normalizing the metric on $\Xx$ once and for all so that the minimal holomorphic
sectional curvature is $-1$, we can summarize the properties
of the Toledo invariant in the following

\begin{prop}\label{prop:tr}
There exists a rational number $\ell_\Xx\in\QQ$ such that 
the Toledo invariant $\tr$ of a representation $\rho:\gg\to G$
has the following properties:
\begin{enumerate}
\item $\tr\in\ell_\Xx\ZZ$;
\item the map $\mathrm{T}:\hom(\gg,G)\to\ell_\Xx\ZZ$ is constant on connected components
of the representation variety, and 
\item $|\tr|\leq|\chi(\sg)|\rkx$, where $\rkx$ is the real rank of $\Xx$.
\end{enumerate}
\end{prop}

\begin{rem} \be
\item The constant $\ell_\Xx$ can be explicitly computed from the
restricted root system of the real Lie group $G$
(see \cite{Burger_Iozzi_Wienhard_kaehler}). In fact, the metric of minimal 
holomorphic sectional curvature $-1$ on $\Xx$ is the $\ell_\Xx$-multiple
of the Bergmann metric given by the bounded domain realization of $\Xx$.
\item The inequality in Proposition~\ref{prop:tr}(3) is due 
to J.~Milnor in the case $G=\sltwo$ \cite{Milnor}
and to V.~Turaev in the case $G=\sp(2n,\RR)$ \cite{Turaev}.
\ee
\end{rem}

We want to give now a very concrete interpretation of the Toledo invariant
of a representation $\gg\to G$ in the case in which $G=\Isom(\Xx)^\circ$,
in terms of generators of $\gg$; this is very much 
in the spirit of Milnor's formula for the Euler number (see \cite{Milnor}).

Let $0\in \Xx$ be a basepoint and $K$ its stabilizer in $G$.
We already alluded to the fact that the center of $K$ is of positive dimension
(see \S~\ref{sec:hss}).
In fact, the $\CC$-vector space structure on the tangent space
$T_0\Xx$ gives an action of $\mathrm{U}(1):=\{z\in\CC:\,|z|=1\}$
which can be ``integrated'', in the sense that there is
a continuous homomorphism 
\bq\label{eq:u_0}
u_0:\mathrm{U}(1)\to K
\eq 
such that for $z\in \mathrm{U}(1)$ 
the differential of the isometry of $\Xx$ defined by 
$u_0(z)$ at $0$ is the multiplication $v\mapsto zv$ for all $v\in T_0\Xx$. 
In particular, since $K$ acts on $T_0\Xx$ faithfully by $\CC$-linear maps, 
the image of $u_0$ is in the center $Z(K)$ of $K$
and in fact coincides with $Z(K)$ when $\Xx$ is irreducible.

Assume hence for the following discussion that $\Xx$ is irreducible.
The homomorphism in (\ref{eq:u_0}) induces a homomorphism 
on the level of fundamental groups 
$\ZZ\to\pi_1(K)=\pi_1(G)$ and hence an isomorphism
\bqn
\ZZ\to\pi_1(G)/\pi_1(G)_{tor}\,.
\eqn
Denoting by $\widehat G$ the covering of $G$ associated to $\pi_1(G)_{tor}$, 
we obtain a topological central extension
\bqn
\xymatrix{
 0\ar[r]
&\ZZ\ar[r]
&\widehat G\ar[r]
&G\ar[r]
&(e)\,.
}
\eqn
The commutator map $\widehat G\times\widehat G\to\widehat G$
factors via $\ZZ$ to give rise to a smooth map
\bqn
\begin{aligned}
G\times G&\to\,\,\,\widehat G\\
(a,b)\,\,&\mapsto \widehat{[a,b]}\,.
\end{aligned}
\eqn
Given a standard presentation of $\gg$,
\bqn
\gg=\bigg\< a_1,b_1,\dots,a_g,b_g:\,\prod_{i=1}^g\,[a_i,b_i]=e\bigg\>\,,
\eqn
to any homomorphism $\rho:\gg\to G$ we can thus associate 
\bqn
\prod_{i=1}^g\,\widehat{[\rho(a_i),\rho(b_i)]}\in\ZZ
\eqn
and the same argument as in \cite{Milnor}, shows that 
\bqn
\tr=\ell_\Xx\prod_{i=1}^g\,\widehat{[\rho(a_i),\rho(b_i)]}\,.
\eqn

As an immediate consequence we have that the Toledo invariant
depends continuously on the representation and hence (1) and
(2) of Proposition~\ref{prop:tr} follow at once. 
The proof of part (3) of the same proposition is more delicate;
one efficient way to prove it uses the value of the simplicial
area of $\sg$ and the value of the norm of the bounded \kahler class
in bounded cohomology determined by Domic and Toledo \cite{Domic_Toledo}
and by Clerc and \O rsted \cite{Clerc_Orsted_2}, 
and will be explained in \S~\ref{sec:tight}.

If $G=\pu(1,1)$ and $\Xx$ is the Poincar\'e disk,
the constant $\ell_\Xx$ in Proposition~\ref{prop:tr} equals 1
and hence for the Toledo invariant of a homomorphism 
\bqn
\rho:\gg\to\pu(1,1)
\eqn
we have that $\tr\in\ZZ$ and $|\tr|\leq2g-2$.
Thus $\tr$ can achieve at most $4g-3$ values.  In fact:

\begin{thm}[Goldman, \cite{Goldman_thesis}, \cite{Goldman_88}]\label{thm:goldman} 
The fibers of
\bqn
\ba
\operatorname{Hom}\big(\gg,\pu(1,1)\big)&\to\{-(2g-2),\dots,0,\dots,2g-2\}\\
\rho\qquad\qquad&\mapsto\qquad\qquad\tr
\ea
\eqn
are exactly the connected components of $\operatorname{Hom}\big(\gg,\pu(1,1)\big)$.
Moreover $\tr=\pm(2g-2)$ if and only if $\rho$ is a hyperbolization,
that is it is faithful with discrete image.
\end{thm}

Let now $G$ be any connected semisimple Lie group with finite center
such that the associated symmetric space $\Xx$ is Hermitian.
Inspired by Proposition~\ref{prop:tr}(3) and by Goldman's result,
we give the following

\begin{defi} We say that a representation $\rho:\gg\to G$
is {\it maximal} if $\tr=\pm\chi(\sg)\rkx$.
\end{defi}

In the sequel we shall show that, in this degree of generality,
maximal representations contain a remarkable amount of structure.
Historically, the following result of H.~Kneser seems to be
the birth certificate of this theme.

\begin{thm}[Kneser, \cite{Kneser}]  Let $f:\sg\to \Sigma_h$ be a continuous map
of compact orientable surfaces, where $h\geq2$.
If $d_f$ is the degree of $f$ we have
\bqn
|d_f(h-1)|\leq|g-1|\,,
\eqn
with equality if and only if $f$ is homotopic to a covering map,
necessarily of degree $d_f$.
\end{thm}

G.~Lusztig observed that this inequality is a consequence of Milnor's inequality
(see \cite{Eells_Wood} for this remark as well as a proof using harmonic maps),
by taking a hyperbolization 
\bqn
\rho:\pi_1(\Sigma_h)\to\pu(1,1)
\eqn
of $\Sigma_h$ and considering the Toledo invariant of the composed
homomorphism 
\bqn
\rho\circ f_\ast:\pi_1(\sg)\to\pu(1,1)\,.
\eqn
Actually, continuing along these lines and applying Goldman's theorem
characterizing maximal representations into $\pu(1,1)$ leads 
to a proof of the equality case in Kneser's theorem.  The reader
might find instructive to fill in the details.

\subsection{Examples of Maximal Representations}
It is a nontrivial and remarkable geometric fact that
any maximal flat in a Hermitian symmetric space can be
``complexified'' hence leading to the existence of maximal polydisks.
\begin{defi} A {\it maximal polydisk} in $\Xx$ is a subHermitian symmetric space 
in $\Xx$ isomorphic to the $r$-fold Cartesian product $(\Dd_{1,1})^r$ of the
Poincar\'e disk, where $r=\rkx$.
\end{defi}

The fact that maximal polydisks exist (and are all conjugate under 
$\Isom(\Xx)^\circ$) \cite[p.~280]{Wolf} is a crucial ingredient in some of the 
examples to follow and their generalization (see Theorem~\ref{thm:ex}).

\begin{exo} The embedding
\bqn
\ba
(\Dd_{1,1})^p\,\,\,&\to\qquad\quad\Dd_{p,q}\\
(z_1,\dots,z_p)&\mapsto\begin{pmatrix}
                           z_1 &      &      &\\
                               &\ddots&      &\\
                               &      &  z_p &\\
                            0  &\hdots&   0  &\\
                         \vdots&      &\vdots&\\
                            0  &\hdots&   0  & 
                       \end{pmatrix}
\ea
\eqn
is isometric and holomorphic and hence defines a maximal polydisk $\Pp\subset\Dd_{p,q}$
associated to the obvious homomorphism 
\bqn
\tau_\Pp:\su(1,1)^p\to\supq\,.
\eqn
\end{exo}

\begin{exo} Let $\RR^{2n}=\RR^2\oplus\dots\oplus\RR^2$ be the direct sum 
of $n$ symplectic planes.  Then the embedding
\bqn
\ba
(\Xx_{\sp(2,\RR)})^n&\to\qquad\Xx_{\sp(2n,\RR)}\\
(J_1,\dots,J_n)&\mapsto\begin{pmatrix}
                       J_1&      &   \\
                          &\ddots&   \\
                          &      &J_n
                       \end{pmatrix}
\ea
\eqn
defines a maximal polydisk $\Pp\subset\Xx_{\sp(2n,\RR)}$
associated to the obvious homomorphism
\bqn
\tau_\Pp:\big(\sp(2,\RR)\big)^n\to\sp(2n,\RR)\,.
\eqn
\end{exo}

\bigskip

Now we present some examples of maximal representations, 
though the proof of their maximality is not necessarily immediate
at this point.

\begin{exo}\label{exo:pd} Let $\Xx$ be any Hermitian symmetric space 
of rank $r$, and $\Pp\subset\Xx$ a maximal polydisk
with associated homomorphism 
\bqn
\tau_\Pp:\su(1,1)^r\to G:=\Isom(\Xx)^\circ\,.
\eqn
Given now $r$ orientation preserving hyperbolizations
$h_1,\dots,h_r:\gg\to\su(1,1)$, the representation
\bqn
\ba
h:\gg&\longrightarrow\quad\su(1,1)^r\\
\gamma&\mapsto\big(h_1(\gamma),\dots,h_r(\gamma)\big)
\ea
\eqn
as well as the composition
\bqn
\tau_\Pp\circ h:\gg\to G
\eqn
are maximal.  
\end{exo}

\begin{exo}\label{exo:tight} Let $h:\gg\to\sltwo$ be an orientation preserving 
hyperbolization and let $\rho_{2n}:\sltwo\to\mathrm{SL}(2n,\RR)$ be the 
$2n$-dimensional irreducible representation of $\sltwo$.  
Since $\rho_{2n}\big(\sltwo\big)$ preserves the standard symplectic 
form on $\RR^{2n}$, we obtain a representation
\bqn
\rho_{2n}\circ h:\gg\to\sp(2n,\RR)
\eqn
which can be proven to be maximal.  
Observe that such a representation $\rho_{2n}\circ h$ 
is in the Hitchin component
$\operatorname{Rep}_H\big(\gg,\sp(2n,\RR)\big)$ (see \S~\ref{sec:intro}).
\end{exo}

The fact that $\tau_\Pp\circ h$ and $\rho_{2n}\circ h$ are maximal depends 
on the property of $\tau_\Pp$ and $\rho_{2n}$ being ``tight homomorphisms'', 
a concept to which we shall return in \S~\ref{sec:tight}.

\begin{exo}  Because of Proposition~\ref{prop:tr}(2), any deformation 
of any of the above representations is maximal.  In particular,
all representations in $\operatorname{Rep}_H\big(\gg,\sp(2n,\RR)\big)$
are maximal.
\end{exo}

\begin{exo}\label{exo:sp(4)} Let $\sg$ be a compact oriented surface of
genus $g\geq2$.  Our objective is to give an example of maximal
representation with Zariski dense image in  $\sp(4,\RR)$
constructed via an explicit deformation of the representation 
in Example~\ref{exo:pd}.
This was triggered by a comment of Toledo, 
who pointed out that a disk diagonally embedded into a polydisk
does not determine uniquely the polydisk.

To this end, write $\sg$ as the sum of two surfaces
$\Sigma_A,\Sigma_B$ 
identified along a separating simple closed loop $\gamma$ on which 
we choose a basepoint $p$
and realize $\gg$ as an amalgamated product
$\Gamma_A\ast_{\<\gamma\>}\Gamma_B$ of $\Gamma_A:=\pi_1(\Sigma_A,p)$ and 
$\Gamma_B:=\pi_1(\Sigma_B,p)$ over the infinite cyclic subgroup $\<\gamma\>$.

 \begin{figure}[!h]
\centering
\psfrag{sigmaa}{$\Sigma_A$}
\psfrag{sigmab}{$\Sigma_A$}
\psfrag{gamma}{$\gamma$}
 \includegraphics[width=1\linewidth]{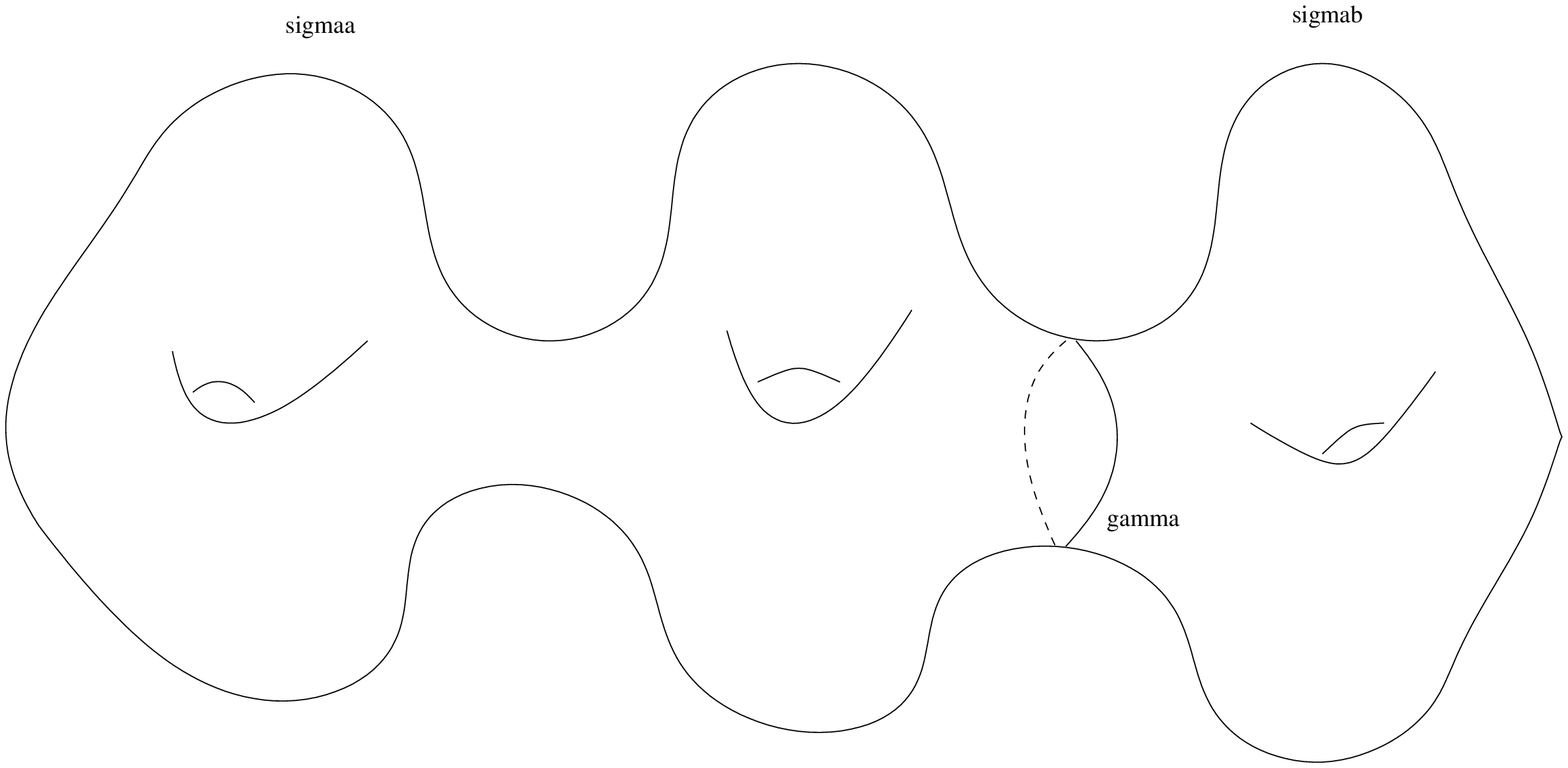}
\caption{\small }
\end{figure}
\vskip1cm

Let $h:\gg\to\mathrm{PSp}(2,\RR)^2$ be defined
by $h(\gamma):=\big(h_1(\gamma),h_2(\gamma)\big)$,
where $h_1,h_2:\gg\to\mathrm{PSp}(2,\RR)$ are two
inequivalent hyperbolizations, and let us choose some lift to
$\sp(2,\RR)^2$, denoted again by $h$ with a small abuse of notation,
\bqn
\ba
h:\gg&\to\quad\sp(2,\RR)^2\\
\gamma&\mapsto\big(h_1(\gamma),h_2(\gamma)\big)\,.
\ea
\eqn
We shall assume that:
\begin{enumerate}
\item[($\ast$)] $h(\gamma)\in\Delta$, where $\Delta$ is the diagonal of $\sp(2,\RR)^2$, and
\item[($\ast\ast$)] the restriction of the two hyperbolic structures 
to $\Sigma_A$ and to $\Sigma_B$ 
are inequivalent,
\end{enumerate}
and we denote by 
\bqn
\begin{aligned}
h_A:\Gamma_A&\to\quad\sp(2,\RR)^2\\
\alpha&\mapsto\big(h_{1,A}(\alpha),h_{2,A}(\alpha)\big)
\end{aligned}
\eqn
and
\bqn
\begin{aligned}
h_B:\Gamma_B&\to\quad\sp(2,\RR)^2\\
\beta&\mapsto\big(h_{1,B}(\beta),h_{2,B}(\beta)\big)
\end{aligned}
\eqn
the restrictions of $h$ to $\Gamma_A$ and $\Gamma_B$ respectively.

Let now $\RR^4=\RR^2\oplus\RR^2$ be the sum of two standard symplectic real planes
(as in ~\ref{subsubsec:sympl}), so that 
\bqn
\sp(4,\RR)=
\left\{g=\begin{pmatrix}A&B\\C&D\end{pmatrix}\in\gl(4,\RR):
{}^tg\begin{pmatrix}J&0\\0&J\end{pmatrix}g=\begin{pmatrix}J&0\\0&J\end{pmatrix}\right\}
\eqn
where $J=\begin{pmatrix}0&1\\-1&0\end{pmatrix}$.
Consider the homomorphism
\bqn
\begin{aligned}
\tau_\Pp:\sp(2,\RR)^2&\to\,\,\,\sp(4,\RR)\\
(A,D)\quad &\mapsto\begin{pmatrix}A&0\\0&D\end{pmatrix}
\end{aligned}
\eqn
associated to the maximal polydisk 
\bqn
\Pp=\left\{\begin{pmatrix}J_1&0\\0&J_2\end{pmatrix}\in\Xx_{\sp(4,\RR)}:\,
J_1,J_2\hbox{ are complex structures on }\RR^2\right\}\,.
\eqn
Then the centralizer $\Zz$ in $\sp(4,\RR)$ of the image $\tau_\Pp(\Delta)$
of the diagonal is 
\bqn
\Zz=\left\{\begin{pmatrix}a\id_2&b\id_2\\c\id_2&d\id_2\end{pmatrix}:
\begin{pmatrix}a&b\\c&d\end{pmatrix}\in\mathrm{O}(2)
\right\}\,.
\eqn
Denoting by $\operatorname{Int}(z)$ the conjugation by $z\in\Zz$, 
the homomorphisms $\tau_\Pp\circ h_A$ and $\operatorname{Int}(z)\circ\tau_\Pp\circ h_B$
coincide on $\<\gamma\>$ (see ($\ast$)).  
Thus by the universal property of amalgams, there is a unique homomorphism
\bqn
\rho_z:\gg\to\sp(4,\RR)
\eqn
whose restriction to $\Gamma_A$ is $\tau_\Pp\circ h_A$ and to $\Gamma_B$ is
$\operatorname{Int}(z)\circ\tau_\Pp\circ h_B$.

\begin{prop}\label{prop:sp(4)} With the above notations:
\begin{enumerate}
\item For every $z\in\Zz^\circ$ the representation $\rho_z:\gg\to\sp(4,\RR)$ 
is maximal, and
\item if $z=\begin{pmatrix}a\id_2&b\id_2\\c\id_2&d\id_2\end{pmatrix}\in\Zz$
satisfies $ac\neq0$, then $\rho_z(\gg)$ is Zariski dense in $\sp(4,\RR)$.
\end{enumerate}
\end{prop}
\begin{proof} (1) If $z=\id$, $\rho_\id$ is the composition of $h:\gg\to\sp(2,\RR)^2$
with $\tau_\Pp:\sp(2,\RR)^2\to\sp(4,\RR)$.  The latter homomorphism
is associated to an embedding realizing $(\Xx_{\sp(2,\RR)})^2$ 
as a maximal polydisk in the Siegel space $\Xx_{\sp(4,\RR)}$ 
and hence $\rho_\id$ is maximal
(see Example~\ref{exo:pd}).
Every $z\in\Zz^\circ$ can be connected to $\id$ by a continuous path,
thus $\rho_z$ is in the same component of $\hom\big(\gg,\sp(4,\RR)\big)$
as $\rho_{\id}$ and thus maximal.

\medskip
\noindent
(2) It follows from ($\ast\ast$) that the image $\rho_z(\gg)$ 
is Zariski dense in the algebraic group $L<\sp(4,\RR)$
generated by $P:=\tau_\Pp\big(\sp(2,\RR)^2\big)$ and $zP z^{-1}$.
Now the condition on $z$ guarantees that the Lie algebra $\frakp$ of $P$ is strictly
contained in the Lie algebra $\frakl$ of $L$.  But it is easily verified
that the representation of $\sp(2,\RR)^2$ on the Lie algebra $\fraks\frakp(4,\RR)$ 
of $\sp(4,\RR)$ obtained by composing $\tau_\Pp$ with the adjoint representation 
is a sum of $\frakp$ and the irreducible representation in dimension 4,
tensor product of the standard $2$-dimensional representation of
$\sp(2,\RR)$ with itself.  This implies that $\frakl=\fraks\frakp(4,\RR)$ 
and hence proves the proposition.
\end{proof}
\end{exo}

\begin{exo} A smooth fiber bundle over a surface $\sg$ with typical fiber
$\Sigma_n$ leads, via the monodromy representation on the first
homology group of the fiber, to a representation
\bqn
\rho:\gg\to\sp(2n,\ZZ)\,.
\eqn
For such representations, D.~Kotschick \cite{Kotschick_sign} showed that
\bqn
|\tr|\leq\frac12|\chi(\sg)|(n-1)\,,
\eqn
and in particular $\rho$ is far from being maximal.
On the other hand D.~Toledo has given examples of
maximal representations into $\sp(4n,\ZZ)$ for all $n\geq1$,
\cite{Toledo_email}.
\end{exo}


\vskip1cm

\section{Tube Type Subdomains and Maximal Representations}\label{sec:res}
Let $G$ be a semisimple real algebraic group with associated symmetric space 
$\Xx$ of Hermitian type.  In view of Goldman's theorem
(see Theorem~\ref{thm:goldman}),
a basic question concerning a maximal representation 
$\rho:\gg\to G$ is whether it is faithful with discrete image.
In addition, when $\Xx$ is not  the Poincar\'e disk,
$\rho(\gg)$ cannot be a lattice in $G$ and thus in this general setting 
there is the question of determining how ``large'' the image of $\rho$
can be.  Concerning the latter question, it is natural to turn
one's attention to the Zariski closure $\lL$ of $\rho(\gg)$. 
While in the preceding examples we have seen that $\lL(\RR)$ can be a 
product $\sp(2,\RR)^r$, or more interestingly $\sp(4,\RR)$,
it turns out however that there are restrictions on $\lL$
and that moreover the determination of these restrictions 
is an essential step in order to answer the question 
about faithfulness and discreteness of $\rho$.

To this purpose, an instructive special case is the family
of Hermitian symmetric spaces of real rank one, 
that is the complex hyperbolic spaces $\Dd_{1,q}$.
This case, beyond $q=1$, was
examined by Toledo and in order to state his result we recall that 
a complex geodesic is an  isomorphic copy of $\Dd_{1,1}$ in
$\Dd_{1,q}$; equivalently, complex geodesics are obtained by taking 
the exponential of a complex line in $T_x\Dd_{1,q}$, for $x\in\Dd_{1,q}$;
they constitute the maximal polydisks in $\Dd_{1,q}$.

\begin{thm}[Toledo, \cite{Toledo_89}]\label{thm:toledo} Any maximal representation
$\rho:\gg\to\mathrm{PU}(1,q)$ stabilizes a complex geodesic.
\end{thm}

Since the stabilizer in $\mathrm{PU}(1,q)$ 
of a complex geodesic is, modulo a compact kernel, 
isomorphic to $\mathrm{PU}(1,1)$, Goldman's theorem \cite{Goldman_thesis} 
applies, and thus $\rho$ is basically obtained via a hyperbolization of $\sg$.

The proof of Toledo's theorem is very much in the spirit
of the Gromov--Thurston proof of Mostow rigidity theorem
and uses notably $\ell^1$-homology and smearing.
(Incidentally, $\ell^1$-homology will play a role
also in our treatment of the Milnor--Wood type inequality in 
Proposition~\ref{prop:tr}(3) as described in \S~\ref{sec:tight}).
A special case of Theorem~\ref{thm:toledo} was already proven  
by Toledo in \cite{Toledo_79} using harmonic mappings techniques.
In the same spirit, 
taking up Hitchin's approach via Higgs bundles,
Bradlow, Garc\'ia-Prada and Gothen made a comprehensive study
of the topology of the connected components of 
$\operatorname{Hom}_{red}\big(\gg,\pu(p,q)\big)$ 
and obtained in particular for the maximal representations
the following

\begin{thm}[Bradlow--Garc\'ia-Prada--Gothen, \cite{Bradlow_GarciaPrada_Gothen}] 
Assume that $1\leq p\leq q$.  If $\rho:\gg\to\pupq$ is maximal 
and reductive, then its image is contained
in $\mathrm{P}\big(\mathrm{U}(p,p)\times\mathrm{U}(q-p)\big)$ up to conjugation.
\end{thm}

This result had been previously obtained by L.~Hern\'andez in the 
case $p=2$, \cite{Hernandez}.
Moreover, an analogous result has been proven 
by Bradlow, Garc\'ia-Prada and Gothen for $\mathrm{SO}^\ast(2n)$, if $n$ is odd
\cite{Bradlow_GarciaPrada_Gothen_conf}.

An equivalent way of stating the theorem asserts
that $\rho(\gg)$ preserves a Hermitian symmetric subspace
of $\Dd_{p,q}$ conjugate to $\Dd_{p,p}$.  
In order to understand the situation in general, 
the relevant concept here is the one of {\it tube type domain}.
For instance, the Hermitian symmetric space associated to 
$\mathrm{PU}(1,1)$ has a realization as upper half plane
but, unlike the bounded domain realization, this type of realization
is not available for all Hermitian symmetric spaces.

\begin{defi} A Hermitian symmetric space $\Xx$ is of {\it tube type} if
it is biholomorphic to a domain
\bqn
\{u+iv:\,u\in V,v\in\Omega\}\subset V\oplus iV
\eqn
where $V$ is a real vector space and $\Omega\subset V$ is a proper open cone.
\end{defi}

\begin{exo}  The space $\Dd_{p,q}$ associated to $\supq$ 
is of tube type if and only if $p=q$, in which case it is biholomorphic to 
\bqn
\operatorname{Herm}_p(\CC)+ i\operatorname{Herm}^+_p(\CC)\,,
\eqn
where $\operatorname{Herm}_p(\CC)$ is the real vector space of Hermitian
matrices and $\operatorname{Herm}^+_p(\CC)$ is the open cone 
of positive definite ones. 
The biholomorphism is given explicitly
by restricting to $\Dd_{p,p}$ the Cayley transformation
\bqn
\ba
\Cc:M_{p,p}(\CC)&\to \,\,\, M_{p,p}(\CC)\\
Z\quad&\mapsto\frac{i(\id+Z)}{(\id-Z)}\,.
\ea
\eqn

The Cayley transformation $\Cc$ sends the real Zariski open subset
of the Shilov boundary $\cs_{p,p}$ of $\Dd_{p,p}$
consisting of matrices $Z\in M_{p,p}(\CC)$ such that 
$\det(I-Z)\neq0$ bijectively into $\operatorname{Herm}_p(\CC)$.
Considering $\operatorname{Herm}^+_p(\CC)$ as an open cone in 
the tangent space of $\operatorname{Herm}_p(\CC)$ at $0$,
and taking its image under the differential at $0$ of $\Cc^{-1}$,
one obtains a cone $\Omega$ in the tangent space $T_{-\id}\cs_{p,p}$
which is invariant under the stabilizer of $-\id$ in $\su(p,p)$.
Translating this cone by the $\su(p,p)$-action, one obtains 
a smooth family of open cones $\Omega_x\subset T_x\cs_{p,p}$,
that is a {\it causal structure} on the Shilov boundary,
which is $\su(p,p)$-invariant.
\end{exo}

\begin{exo} Under the same Cayley transformation, 
the bounded domain realization of the Siegel space is sent biholomorphically to 
\bqn
\operatorname{Sym}_n(\RR)+ i\operatorname{Sym}_n^+(\RR)\,,
\eqn
where $\operatorname{Sym}_n(\RR)$ denotes the vector space consisting
of real symmetric $n\times n$ matrices, and realizes
the symmetric space associated to $\sp(2n,\RR)$ as a tube type domain.

The open cone $\operatorname{Sym}^+(\RR)$ defines in the same way 
as above an $\sp(V)$-invariant causal structure on $\Ll(V)$
which will be used in the proof of Corollary~\ref{cor:rect}.
\end{exo}

The examples above serve to illustrate the general fact 
that the Shilov boundary of a bounded symmetric domain of 
tube type admits an invariant causal structure.  
Among bounded symmetric domains, this property characterizes those of tube type.
The general classification of Hermitian symmetric spaces 
relative to the notion of tube type is as follows:
\bigskip
{
\begin{center}
\noindent
\begin{tabular}{c|c}
\hline
{\it Tube Type} & {\it Non-Tube Type}\\
\hline
$\su(p,p)$&$\su(p,q)$\\
&$p<q$\\
\hline
$\sp(2n,\RR)$& \\
\hline
$\mathrm{SO}^\ast(2n)$&$\mathrm{SO}^\ast(2n)$\\
$n$ even& $n$ odd\\
\hline
$\mathrm{SO}(2,n)$& \\
\hline
$E_7(-25)$&$E_6(-14)$\\
\hline
\end{tabular}
\end{center}
}
\bigskip
\noindent
where $E_7(-25)$ and $E_6(-14)$ correspond to the exceptional Hermitian symmetric spaces
of complex dimension $27$ and $16$ respectively. 

An essential feature of a Hermitian symmetric space of rank $\rkx$ is that
(holomorphically embedded) maximal Hermitian subdomains of tube type always exist, 
are of rank $\rkx$ and are all conjugate.

Notice that in the rank one case, that is for complex hyperbolic $n$-space,
the notion of maximal tube type subdomain and maximal polydisk coincide.
This ambiguity left open the correct generalization of Toledo's theorem
until the construction of a maximal representation with Zariski dense image
in a tube type domain \cite{Burger_Iozzi_Wienhard_ann}, 
of which Example~\ref{exo:sp(4)} is a particular case.

We can finally state the structure theorem for maximal representations. 
\begin{thm}[\cite{Burger_Iozzi_Wienhard_ann}, \cite{Burger_Iozzi_Wienhard_tol}]\label{thm:main}
Let $\gG$ be a connected semisimple real algebraic group and assume that 
the symmetric space $\Xx_G$ associated to $G:=\gG(\RR)^\circ$ is Hermitian
symmetric.
Let $\rho:\gg\to G$ be a maximal representation.  Then:
\begin{enumerate}
\item $\gg$ acts properly discontinuously, via $\rho$, on $\Xx_G$;
\item the Zariski closure $\lL$ of $\rho(\gg)$ is reductive;
\item the connected component $L:=\lL(\RR)^\circ$ stabilizes a maximal tube type
 subdomain $T\subset\Xx_G$;
\item the symmetric space $\Xx_L$ associated to $L$ is Hermitian
 of tube type and the isometric embedding $\Xx_L\hookrightarrow T$ is tight.
\end{enumerate}
\end{thm}

We illustrate the above theorem in the examples of \S~\ref{sec:inv}.
For the notion of {\it tight embedding} see \S~\ref{sec:tight}.

\begin{exo*}\ref{exo:pd} The orientation preserving hyperbolizations 
$h_1,\dots,h_r$ fall into $s$ equivalence classes
modulo $\su(1,1)$-equivalence, with $1\leq s\leq r$;  
the group $L$ is then isomorphic
to $\su(1,1)^s$ and $\Xx_L$ is a product of $s$-subdiagonals in the
maximal polydisk $\Pp$.  Any diagonal disk in $\Pp$ determines 
the same maximal tube type domain $T$ and $\Xx_L\subset\Pp\subset T$.
The first statement of the theorem is obvious in this example. 
\end{exo*}

\begin{exo*}\ref{exo:tight} In this case $L\cong\sp(2,\RR)$, $T=\Xx_{\sp(2n,\RR)}$,
and $\Xx_L\subset\Xx_{\sp(2n,\RR)}$ is a geodesically embedded disk,
holomorphic if and only if $n=1$.  Item (1) of the theorem is obvious
also in this case.
\end{exo*}

\begin{exo*}\ref{exo:sp(4)} Under the conditions of Proposition~\ref{prop:sp(4)}(2),
the maximal representation
\bqn
\rho_z:\gg\to\sp(4,\RR)
\eqn
has Zariski dense image, so $L=\sp(4,\RR)$ and 
$\Xx_L=T=\Xx_{\sp(4,\RR)}$ by construction.  
On the other hand, Theorem~\ref{thm:main}(1) implies that 
$\rho_z$ is injective with discrete image, 
a fact that in this case is not at all obvious from the construction.
\end{exo*}

As alluded to earlier, Example~\ref{exo:sp(4)} is a particular case
of a general fact which we now state.  Let $r=\rkx$, 
let $\Pp\subset\Xx$ be a maximal polydisk and
\bqn
\tau_\Pp:\su(1,1)^r\to\Isom(\Xx)^\circ\,.
\eqn
the associated homomorphism.

\begin{thm}[\cite{Burger_Iozzi_Wienhard_ann}, \cite{Burger_Iozzi_Wienhard_tol}]\label{thm:ex} Assume that $\Xx$ is of tube type 
and let $\rho_0:\gg\to\Isom(\Xx)^\circ$
be the maximal representation obtained by composing
a hyperbolization of $\gg\to\su(1,1)$ with the diagonal
embedding of $\su(1,1)$ in $\su(1,1)^r$ followed by $\tau_\Pp$.
Then $\rho_0$ admits a continuous deformation
$\rho_t:\gg\to\Isom(\Xx)^\circ$, for $t\geq0$, such that 
$\rho_t(\gg)$ is Zariski dense in $\Isom(\Xx)^\circ$ for $t>0$.
\end{thm}

Observe that $\rho_t$, being a continuous deformation of a maximal 
representation, is maximal as well.


\vskip1cm

\section{Tight Homomorphisms}\label{sec:tight}
A fundamental role in the study of maximal representations
of surface groups is played by tight homomorphisms,
which generalize maximal representations of surface groups,
in that it is a notion defined for any continuous homomorphism 
of a locally compact group into the group of isometries of
a Hermitian symmetric space.

The definition of tight homomorphism rests on basic concepts
in bounded continuous cohomology which we briefly recall;
for a comprehensive treatment see \cite{Monod_book} and \cite{Burger_Monod_GAFA}.
We start with the more familiar concept of continuous group cohomology.
For a locally compact group $G$, its continuous cohomology 
$\hc^\bu(G,\RR)$ is the cohomology of the complex 
$\big(\mathrm{C}(G^\bu,\RR)^G,d^\bu\big)$ of $G$-invariant 
real valued continuous cochains,
where $d^\bu$ is the usual homogeneous coboundary.
The bounded continuous cohomology $\hcb^\bu(G,\RR)$
is then the cohomology of the subcomplex $\big(\cb(G^\bu,\RR)^G,d^\bu\big)$
of $G$-invariant bounded continuous cochains.  
The complex $\big(\cb(G^\bu,\RR)^G,d^\bu\big)$ equipped with the supremum
norm is a complex of Banach spaces with continuous coboundary operators,
and hence $\hcb^\bu(G,\RR)$ is endowed with a quotient seminorm.
Also, the inclusion of the complex of bounded continuous functions
into the one of continuous functions gives rise to a comparison map
\bq\label{eq:comp}
{\tt c}^\bu_G:\hcb^\bu(G,\RR)\to\hc^\bu(G,\RR)
\eq
which encodes subtle properties of $G$ of geometric and algebraic nature.
See \cite{Bavard}, \cite{Ghys_87}, \cite{Mineyev}, \cite{Burger_Monod_JEMS},
\cite[\S~V.13]{Burger_Monod_GAFA}, \cite{Burger_Iozzi_supq}, and 
also \cite{Brooks}, \cite{Brooks_Series}, \cite{Grigorchuk}, 
\cite{Soma}, \cite{Bestvina_Fujiwara}, \cite{Epstein_Fujiwara}, 
\cite{Gambaudo_Ghys}, \cite{Entov_Polterovich}, \cite{Biran_Entov_Polterovich}, 
\cite{Kotschick_not}, \cite{Kotschick_proc}
in relation with the existence of quasi-morphisms.

If now $G$ is a connected semisimple Lie group with finite center
and associated symmetric space $\Xx$, we have seen that the
complex $\big(\Omega^\bu(\Xx)^G,d^\bu\big)$ 
of $G$-invariant differential forms on $\Xx$
coincides with its cohomology (see Lemma~\ref{lem:closed}) 
and, in fact, there is a canonical isomorphism \cite{Van_Est}
\bqn
\hc^\bu(G,\RR)\cong\Omega^\bu(\Xx)^G\,.
\eqn

Let us now specialize to the case of interest to us, namely
when $\Xx$ is Hermitian symmetric and $\omega_\Xx\in\Omega^2(\Xx)^G$
is its \kahler form.  A continuous cocycle defining the class 
$\kx\in\hc^2(G,\RR)$ corresponding to $\ox$ is
\bqn
c_\Xx(g_1,g_2,g_3)=\int_{\Delta(g_1x,g_2x,g_3x)}\ox\,,
\eqn
where $x\in \Xx$ is a basepoint and $\Delta(g_1x,g_2x,g_3x)$ is any 
smooth two-simplex with geodesic sides and vertices $g_1x,g_2x,g_3x$. 

There is a general conjecture of Dupont to the extent that cocycles 
obtained by integrating $G$-invariant differential forms (of any degree) 
should be bounded, \cite{Dupont}. 
In terms of the comparison map, this suggests the following

\begin{que*} Let $G$ be a connected semisimple Lie group with finite center.
Is the comparison map (\ref{eq:comp}) surjective in all degrees?
\end{que*}

This turns out to be true for forms representing specific classes 
(see \cite{Gromov_82}, \cite{Dupont}, \cite{Savage}, \cite{Bucher_thesis},
\cite{Lafont_Schmidt}) and in particular for the \kahler form
was first shown by Dupont \cite{Dupont}.  
In fact, with the assumed normalization on the metric of $\Xx$
(see \S~\ref{sec:inv}), one has the equality
\bq\label{eq:bound}
\|c_\Xx\|_\infty=\pi\rkx
\eq
due to Domic and Toledo for classical domains \cite{Domic_Toledo}
and to Clerc and \O rsted in the general case \cite{Clerc_Orsted_2}.
Thus $c_\Xx$ defines a continuous bounded class $\kxb\in\hcb^2(G,\RR)$
to which we shall refer to the {\it bounded \kahler class};
and for which one has the following theorem (see also Proposition~\ref{prop:maslov}):

\begin{thm}[Domic--Toledo \cite{Domic_Toledo}, Clerc--\O rsted \cite{Clerc_Orsted_2}]\label{thm:norm}
If the metric on $\Xx$ is normalized to have minimal holomorphic sectional curvature 
$-1$, then 
\bqn
\|\kxb\|=\pi\rkx\,.
\eqn
\end{thm}

Given now locally compact groups $H$ and $G$, any continuous homomorphism
$\rho:H\to G$ induces canonical pullbacks $\rho^\bu$ and 
$\rho^\bu_{\mathrm{b}}$ respectively in continuous and bounded
continuous cohomology, by precomposition of continuous (bounded)
cochains with $\rho$; the resulting linear maps have the property
that the diagram
\bq\label{eq:diagram}
\xymatrix{
 \hcb^\bu(G,\RR)\ar[r]^{\rho^\bu_{\mathrm{b}}}\ar[d]_{{\tt c}_G^\bu}
&\hcb^\bu(H,\RR)\ar[d]^{{\tt c}_H^\bu}\\
 \hc^\bu(G,\RR)\ar[r]^{\rho^\bu}
&\hc^\bu(H,\RR)
}
\eq
commutes, and moreover the pullback in bounded continuous cohomology is
norm decreasing, namely for all $\alpha\in\hcb^n(G,\RR)$,
\bqn
\|\rho^{(n)}_{\mathrm{b}}(\alpha)\|\leq\|\alpha\|\,.
\eqn

\begin{defi}[\cite{Burger_Iozzi_Wienhard_tight}, \cite{Wienhard_thesis}] 
Let $G$ be a connected semisimple group with finite center
and such that the associated symmetric space $\Xx$ is Hermitian,
and let $H$ be any locally compact group.
A continuous homomorphism $\rho:H\to G$ is {\it tight}
if it preserves the norm of the bounded \kahler class, that is if
\bqn
\big\|\rho^{(2)}_\mathrm{b}(\kxb)\big\|=\big\|\kxb\big\|\,.
\eqn
\end{defi}

To motivate this definition, we sketch a proof 
of the inequality in Proposition~\ref{prop:tr}(3).
Since $\sg$ is a $K(\gg,1)$, we have in particular
a canonical isomorphism 
\bqn
\h^2(\gg,\RR)\to\h^2(\sg,\RR)
\eqn
which, if $\kx\in\hc^2(G,\RR)$ allows us to see $\rho^{(2)}(\kx)\in\h^2(\gg,\RR)$
as a singular class in $\h^2(\sg,\RR)$ and evaluate it on the fundamental class
$[\sg]\in\h_2(\sg,\RR)$ of $\sg$; recall that $\sg$ is oriented
once and for all.  Then if 
\bqn
\<\,\cdot\,,\,\cdot\,\>:
\h^2(\sg,\RR)\times\h_2(\sg,\RR)\to\RR
\eqn
denotes the pairing, analogously to the classical case of 
the Euler number, we have
\bqn
\tr=\frac{1}{2\pi}\big\<\rho^{(2)}(\kx),[\sg]\big\>\,.
\eqn
The proof of the Milnor--Wood type inequality in Proposition~\ref{prop:tr}(3)
will follow from the interpretation of this invariant in bounded cohomology.
To this purpose, following Gromov \cite{Gromov_82}, recall that 
the $\ell^1$-homology of $\sg$ is the homology $\h_{\bu,\ell^1}(\sg,\RR)$
of the complex of singular $\ell^1$-chains, while the bounded cohomology
$\hb^\bu(\sg,\RR)$ is the cohomology of the dual Banach space complex; 
consequently, $\ell^1$-homology and bounded cohomology acquire quotient
seminorms and there is the canonical pairing
\bqn
\<\,\cdot\,,\,\cdot\,\>_\mathrm{b}:
\hb^n(\sg,\RR)\times\h_{n,\ell^1}(\sg,\RR)\to\RR
\eqn
which satisfies the property that for all $\alpha\in\hb^n(\sg,\RR)$
and all $a\in\h_{n,\ell^1}(\sg,\RR)$
\bqn
|\<\alpha,a\>_{\mathrm{b}}|\leq\|\alpha\|\,\|a\|_{\ell^1}\,.
\eqn
These notions have been introduced by Gromov for any topological space $X$
and one has the Gromov--Brooks canonical isometric isomorphism
(see \cite{Gromov_82} and \cite{Brooks})
\bqn
\hb^\bu\big(\pi_1(X),\RR\big)\cong\hb^\bu(X,\RR)\,,
\eqn
a rather deep fact depending on higher homotopy groups being
Abelian and hence amenable.  
In our situation one can explicitly write an isometric isomorphism
\bq\label{eq:gromov_brooks}
\hb^2(\gg,\RR)\cong\hb^2(\sg,\RR)
\eq
compatible with the isomorphism in ordinary cohomology,
by choosing a hyperbolic metric on $\sg$ and using the technique
of straightening simplices.

Starting now with the bounded \kahler class $\kxb\in\hcb^2(G,\RR)$,
and applying the pullback in bounded cohomology and the isomorphism
(\ref{eq:gromov_brooks}), we obtain the class 
$\rho^{(2)}_\mathrm{b}(\kxb)\in\hb^2(\sg,\RR)$
which corresponds, using (\ref{eq:diagram}), to $\rho^{(2)}(\kx)\in\h^2(\sg,\RR)$
under the comparison map in singular cohomology
\bqn
\hb^2(\sg,\RR)\to\h^2(\sg,\RR)\,.
\eqn
This latter being the dual of the natural map
\bq\label{eq:comp_hom}
\h_2(\sg,\RR)\to\h_{2,\ell^1}(\sg,\RR),
\eq
we have that
\bqn
\big\<\rho^{(2)}(\kx),[\sg]\big\>=
\big\<\rho^{(2)}_{\mathrm{b}}(\kxb),[\sg]_{\ell^1}\big\>_{\mathrm{b}}\,,
\eqn
where $[\sg]_{\ell^1}$ denotes the image of $[\sg]$ under (\ref{eq:comp_hom}).
Thus 
\bqn
|\tr|\leq\frac{1}{2\pi}\big\|\rho^{(2)}_\mathrm{b}(\kxb)\big\|\,\big\|[\sg]\big\|_{\ell^1}\,.
\eqn
Recall now that the $\ell^1$-norm $\big\|[\sg]\big\|_{\ell^1}$ of the fundamental
class is called the simplicial area of $\sg$ and, by \cite{Gromov_82},
\bqn
\big\|[\sg]\big\|_{\ell^1}=4g-4\,.
\eqn
This, together with the norm decreasing property of the pullback
in bounded cohomology and the value of the norm of the \kahler class
in Theorem~\ref{thm:norm},
implies on the one hand the inequality in Proposition~\ref{prop:tr}(3)
and on the other the following

\begin{prop}
Any maximal representation is a tight homomorphism.
\end{prop}

The following general result about tight homomorphisms, together
with the above proposition, implies part of Theorem~\ref{thm:main}.

\begin{thm}[\cite{Burger_Iozzi_Wienhard_tight}, \cite{Wienhard_thesis}]\label{thm:tight}  Let $\gG$ be a semisimple real algebraic group,
$H$ a locally compact group, and assume that the symmetric space
$\Xx_G$ associated to $G:=\gG(\RR)^\circ$ is Hermitian.  Then 
for a tight homomorphism $\rho:H\to G$ the following holds:
\be
\item The Zariski closure $\lL$ of the image $\rho(H)$ 
is reductive;
\item the real reductive group $L:=\lL(\RR)^\circ$ has compact 
centralizer in $G$; and
\item the symmetric space $\Xx_L\subset\Xx_G$ associated to $L$ is Hermitian.
\ee
\end{thm}

Notice that the totally geodesic embedding $\Xx_L\subset\Xx_G$ in (3) 
is not necessarily holomorphic.  However, there is a notion of tight embedding for 
Hermitian symmetric spaces which parallels the one for homomorphisms.

\begin{defi}[\cite{Burger_Iozzi_Wienhard_tight}, \cite{Wienhard_thesis}] 
Given a totally geodesic embedding
\bqn
f:\Yy\to\Xx
\eqn
of Hermitian symmetric spaces, we say that $f$ is {\it tight}
if
\bqn
 \sup_{\Delta\subset\Yy}\int_\Delta f^\ast\ox
=\sup_{\Delta\subset\Xx}\int_\Delta \ox
\eqn
where the supremum is taken over all smooth triangles with geodesic sides.
\end{defi}

This corresponds for the associated homomorphism 
\bq\label{eq:tight-homo}
\rho:H_\Yy\to\Isom(\Xx)^\circ\,,
\eq
where $H_\Yy$ is an appropriate finite covering of $\Isom(\Yy)^\circ$,
to be a tight homomorphism.  With this terminology, the inclusion
$\Xx_L\subset\Xx_G$ in Theorem~\ref{thm:tight}(3) is a tight embedding.

Here are some examples of tight embeddings:

\begin{exo}\label{exo:1}  The homomorphism $\tau_\Pp:\su(1,1)^r\to\Xx$
associated to a maximal polydisc $\Pp\subset\Xx$ is tight;
evidently, the embedding $\Pp\subset\Xx$ is both holomorphic and tight.
\end{exo}

\begin{exo}\label{exo:2} The irreducible representation 
$\rho_{2n}:\sltwo\to\sp(2n,\RR)$  is tight, and the associated
totally geodesic embedding of $\Hh_\RR^2\to\Xx_{\sp(2n,\RR)}$
is a tight embedding which is holomorphic only if $n=1$.
\end{exo}

\begin{exo} The embedding $T\subset\Xx$ of a maximal tube type subdomain
in $\Xx$ is tight and holomorphic.
\end{exo}

\begin{exo} If $\Yy$ is an irreducible Hermitian symmetric space 
and $f:\Yy\to\Xx$ is a totally geodesic embedding, then $f$
is tight if and only if 
\bqn
f^\ast\ox=\pm\frac{\rkx}{\rky}\omega_\Yy\,.
\eqn
\end{exo}

\begin{exo} The embedding $\Xx_{\sp(V)}\to\Xx_{\su(V_\CC)}$
in Example ~\ref{subsubsec:sympl} is tight and holomorphic.
\end{exo}

Observe now the following simple

\begin{prop} Let $H,G$ be connected semisimple Lie groups
with finite center and associated symmetric spaces of Hermitian type.
If $\rho:\gg\to H$ is maximal and $\rho':H\to G$ is tight,
then $\rho'\circ\rho$ is maximal.
\end{prop}

This, together with Examples~\ref{exo:1} and \ref{exo:2} above
justifies the maximality of the representations in 
Examples~\ref{exo:pd} and \ref{exo:tight}.

Notice that in general totally geodesic embeddings between bounded
symmetric domains do not induce maps between the corresponding
Shilov boundaries even if they are holomorphic.  
This is however something else that tight homomorphism can provide for us, namely

\begin{thm}[\cite{Burger_Iozzi_Wienhard_tight}, \cite{Wienhard_thesis}]\label{thm:shilov} Let $\Xx,\Yy$ be Hermitian symmetric spaces and
$f:\Yy\to\Xx$ a tight embedding with associated homomorphism
$\rho:H_\Yy\to\Isom(\Xx)^\circ$ (see (\ref{eq:tight-homo})).
Then there exists a $\rho$-equivariant map
\bqn
\check f:\cs_\Yy\to\cs_\Xx\,.
\eqn
\end{thm}

Remark that, since the Shilov boundary of a Hermitian symmetric space
is a homogeneous space, if such $\rho$-equivariant map exists, it is
unique (up to translations).

The above theorem allows us also to deduce in great generality 
the existence of boundary maps for tight homomorphisms.  
Let $\Lambda$ be a countable discrete group and 
let $\theta$ be a probability measure on $\Lambda$.
Recall that a {\it Poisson boundary} of the pair $(\Lambda,\theta)$
is a measurable $\Lambda$-space $B$ with a quasiinvariant probability measure $\nu$ 
such that there exists an isometric isomorphism between the space of bounded
$\theta$-harmonic functions
\bqn
\ba
\Hh^\infty(\Lambda,\theta):=\big\{f:\Lambda\to\RR:\,&f\hbox{ is bounded and }\\
 &f(g)=\int_\Lambda f(gh)d\theta(h),\forall g\in\Lambda\big\}
\ea
\eqn
and the space $\linfty(B,\nu)$, given by the Poisson formula
\bq\label{eq:poisson}
f(g)=\int_B\psi(gx)d\nu(x)\,.
\eq
Although we shall not need it here, we recall that, under natural
assumptions on the measure $\theta$, a Poisson boundary in fact exists
even for locally compact second countable groups, 
\cite{Kaimanovich_poisson}. 

An immediate consequence of the Poisson formula (\ref{eq:poisson})
is that the measure $\nu$ is $\theta$-stationary, that is $\theta\ast\nu=\nu$.
Moreover, it will be essential for our purposes that the action of $\Lambda$ 
on the Poisson boundary $B$ is amenable 
with respect to the measure $\nu$, \cite{Zimmer:Poisson}.

\begin{thm}\label{thm:map}
Let $\Lambda$ be a countable discrete group with probability measure $\theta$
and let $\gG$ be a semisimple real algebraic group such that 
the symmetric space $\Xx$ associated to $G:=\gG(\RR)^\circ$ is Hermitian.
If $(B,\nu)$ is a Poisson boundary for $(\Lambda,\theta)$
and $\rho:\Lambda\to G$ is a tight homomorphism,
then there exists a $\rho$-equivariant measurable map 
\bqn
\varphi:B\to\cs_\Xx\,.
\eqn
\end{thm}

\begin{proof} Let $\lL$ be the Zariski closure of $\rho(\Lambda)$.
By Theorem~\ref{thm:tight} the symmetric space $\Yy$ associated to 
$L:=\lL(\RR)^\circ$ is Hermitian symmetric and the embedding $\Yy\to\Xx$ is tight, 
so that Theorem~\ref{thm:shilov} implies the existence of a $\rho$-equivariant 
map $\check f$ between the corresponding Shilov boundaries
\bq\label{eq:shilov}
\check f:\cs_\Yy\to\cs_\Xx\,.
\eq
Let $\qQ<\lL$ be a maximal parabolic subgroup defined over $\RR$ 
such that $\cs_\Yy\cong\lL(\RR)/\qQ(\RR)$, 
and let $\pP<\qQ$ be a minimal parabolic subgroup
defined over $\RR$ contained in $\qQ$, 
so that we have an equivariant map 
\bq\label{eq:surj-par}
\lL(\RR)/\pP(\RR)\twoheadrightarrow\lL(\RR)/\qQ(\RR)\cong\cs_\Yy\,.
\eq

Since the action of $\Lambda$ on $(B,\nu)$ is amenable, 
there exists a $\rho$-equivariant measurable map 
\bqn
\varphi_0:B\to\Mm^1\big(\lL(\RR)/\pP(\RR)\big)
\eqn
where $\Mm^1\big(\lL(\RR)/\pP(\RR)\big)$ denotes 
the space of probability measures on $\lL(\RR)/\pP(\RR)$.
Since $\rho:\Lambda\to L$ has Zariski dense image,
the $\Lambda$-action on $\lL(\RR)/\pP(\RR)$ is mean proximal
(see \cite[Theorem~7.3]{Burger_Iozzi_supq}): 
this, together with the fact that $\nu$ is $\theta$-stationary,
implies that for $\nu$-a.\,e. $b\in B$, $\varphi_0(b)$ is a Dirac measure,
thus providing a map
\bqn
\varphi_0:B\to\lL(\RR)/\pP(\RR)
\eqn
which composed with the maps in (\ref{eq:shilov}) and (\ref{eq:surj-par})
provides the required $\rho$-equivariant map.
\end{proof}

As an application, given a compact surface group $\gg$,
choose a hyperbolization of $\sg$ and let $\Gamma$ be the realization 
of $\gg$ as a cocompact lattice in $\pu(1,1)$.
Then $\Gamma$ acts naturally on $\SS^1=\partial\Dd_{1,1}$ and, 
in fact, a theorem of Furstenberg asserts that there exists
a probability measure $\theta$ on $\Gamma$ such that $\SS^1$ with the
Lebesgue measure $\lambda$ is a Poisson boundary of $(\Gamma,\theta)$.

\begin{cor}\label{cor:boundary-map}  Let $\gG$ be a semisimple 
real algebraic group such that the symmetric space
$\Xx$ associated to $G:=\gG(\RR)^\circ$ is Hermitian
and let $\rho:\Gamma\to G$ be a tight embedding of a cocompact 
lattice $\Gamma<\mathrm{PU}(1,1)$.  Then there exists a $\rho$-equivariant
measurable map 
\bqn
\varphi:\SS^1\to\cs_\Xx\,.
\eqn
\end{cor}

\begin{rem}\label{rem:way-out} 
For technical purposes one can show that if $F\subset\cs_\Xx$
is the set of points which are not transverse to a given point in $\cs_\Xx$,
and $\varphi$ is the map in Corollary~\ref{cor:boundary-map}, 
the set $\varphi^{-1}(F)$ has Lebesgue measure zero in $\SS^1$.
\end{rem}

\vskip1cm

\section{Symplectic Anosov Structures}\label{sec:anosov}
We focus in this section on maximal representations 
into a symplectic group $\sp(V)$. Let thus
$\rho:\gg\to\sp(V)$ be any representation. 
We choose a hyperbolization $\Sigma$ of $\sg$, 
and let $\Gamma<\pu(1,1)=\operatorname{Aut}(\Dd_{1,1})^\circ$
be the resulting realization of $\gg$.
From now on we consider $\rho$ as a representation of $\Gamma$.
The geodesic flow $\tilde{g_t}$ on the unit tangent bundle
$T^1\Dd_{1,1}$ gives rise to a flow $\tilde{g_t}^\rho$ on the total
space of the flat symplectic bundle 
$\widetilde E^\rho:=T^1\Dd_{1,1}\times V$ over $T^1\Dd_{1,1}$ commuting 
with the diagonal $\Gamma$-action 
given by $\gamma(u,x):=\big(\gamma u,\rho(\gamma)x\big)$
which hence descends to a flow $g_t^\rho$ on the quotient
$E^\rho:=\Gamma\backslash(T^1\Dd_{1,1}\times V)$ which is a flat symplectic bundle 
over the unit tangent bundle $T^1\Sigma$.  The projection
\bq\label{eq:p}
p:E^\rho\to T^1\Sigma
\eq
is then equivariant with respect to the $g_t^\rho$-action
on $E^\rho$ and to the action of the geodesic flow $g_t$ on $T^1\Sigma$.

Let $\<\,\cdot\,,\,\cdot\,\>:E^\rho\times_p E^\rho\to\RR$
be the symplectic form on $E^\rho$.  A {\it positive complex structure} on the
symplectic bundle is a continuous section 
\bqn 
J:T^1\Sigma\to\operatorname{End}(E^\rho)
\eqn
such that 
\be
\item $J_u$ is a complex structure on the fiber $E^\rho(u)$, and
\item the form $\<\,\cdot\,,J\cdot\,\>$ is symmetric and positive
definite in each fiber.
\ee
We denote by $\|\cdot\|:E^\rho\to\RR_+$ the resulting Euclidean norm,
and by $\|\cdot\|_u$ its value on the fiber $E^\rho(u)$ above the
point $u\in T^1\Sigma$.

Observe that any symplectic bundle over a paracompact base admits 
a positive complex structure.  A {\it Lagrangian subbundle} of a symplectic bundle
is a subbundle such that each fiber is a Lagrangian subspace.
With this terminology we have then the following
\begin{thm}\label{thm:anosov}
Assume that $\rho:\Gamma\to\sp(V)$ is a maximal representation. 
Then there is a $g_t^\rho$-invariant splitting
\bqn
E^\rho=E_-^\rho\oplus E_+^\rho
\eqn
into continuous Lagrangian subbundles, and there exist
a positive complex structure $J$ and a constant $A>0$ such that
\be
\item $J$ interchanges $E_-^\rho$ and $E_+^\rho$, and
\item for all $t\geq0$, 
\bqn 
\|g_t^\rho\xi\|\leq e^{-At}\|\xi\|\hfill\hbox{ for all }\xi\in E_+^\rho
\eqn
and
\bqn
\|g_{-t}^\rho\xi\|\leq e^{-At}\|\xi\|\hfill\hbox{ for all }\xi\in E_-^\rho\,.
\eqn
\ee
\end{thm}
This result has interesting consequences on the metric properties
of a maximal representation.  To describe them, as well as for
convenience in the proofs in \S~\ref{sec:anproofs}, 
we specify a left invariant metric 
on the symmetric space $\Xx_{\sp(V)}$ associated to $\sp(V)$.  
Recall that $\Xx_{\sp(V)}$ is the set of complex structures $J$ on $V$
such that $\<\,\cdot\,,J\cdot\,\>$ is symmetric and positive definite.
Denoting by $q_J$ the corresponding Euclidean norm on $V$,
and by $\|\id\|_{J_1,J_2}$ the norm of the identity map
between $(V,q_{J_1})$ and $(V,q_{J_2})$, we set 
\bqn
d(J_1,J_2):=\big|\ln\|\id\|_{J_1,J_2}\big|+\big|\ln\|\id\|_{J_2,J_1}\big|
\qquad J_1,J_2\in\Xx_{\sp(V)}\,.
\eqn
Of course, this distance is equivalent to the $G$-invariant 
Riemannian distance on $\Xx_{\sp(V)}$, but it is more convenient for our purposes.

The statement of next corollary does not depend on the choice of a 
hyperbolization.

\begin{cor}\label{cor:qi} Let $\rho:\gg\to\sp(V)$ be a maximal representation,
$J\in\Xx_{\sp(V)}$ a basepoint and $\ell$ the word length on $\gg$.
Then the orbit map
\bqn
\ba
\rho_J:\gg&\to\Xx_{\sp(V)}\\
\gamma\,&\mapsto\rho(\gamma) J
\ea
\eqn
is a quasiisometric embedding, that is there are constants $A,B>0$ 
such that for every $\gamma\in\Gamma$
\bqn
A^{-1}\ell(\gamma)-B\leq d\big(\rho(\gamma)J,J\big)\leq A\ell(\gamma)+B\,.
\eqn
\end{cor}

\bigskip

Essential in the proof of Theorem~\ref{thm:anosov} is the existence 
of the boundary map obtained in Corollary~\ref{cor:boundary-map} 
from the boundary $\SS^1=\partial\Dd_{1,1}$ of the Poincar\'e disk 
into the space of Lagrangians $\Ll(V)$ 
which relates the Maslov cocycle (see \S~\ref{sec:proofs}) to the orientation cocycle
on $\SS^1$. 
{\sl A priori} this map is only measurable, but as a consequence of the
continuity of the splitting in Theorem~\ref{thm:anosov},
it turns out to be continuous.  In fact, this map plays a role
analogous to the one of hyperconvex curves in the study of 
the Hitchin component of $\operatorname{Hom}\big(\gg,\mathrm{SL}(n,\RR)\big)$
in \cite{Labourie_anosov}.

\begin{cor}\label{cor:rect} Let $\rho:\Gamma\to\sp(V)$ be 
a maximal representation.  Then there is a $\rho$-equivariant
continuous injective map 
\bqn
\varphi:\SS^1\to\Ll(V)
\eqn
with rectifiable image.
\end{cor}


\vskip1cm

\section{Bounded Cohomology at Use}\label{sec:proofs}

The definition of continuous bounded cohomology in \S~\ref{sec:tight}
is not very useful from a practical point of view,
as many natural cocycles of geometric origin are not continuous.
The homological algebra approach developed in \cite{Burger_Monod_GAFA},
\cite{Monod_book}, \cite{Burger_Iozzi_app} and \cite{Burger_Iozzi_def}
allows us to overcome this obstacles in the usual way:
as in the homological algebra approach to continuous cohomology,
there are appropriate notions of coefficients modules, 
of relatively injective modules and of strong resolutions,
that is resolutions with an appropriate homotopy operator.
The underlying philosophy is that we need not restrict to the
standard resolution in \S~\ref{sec:tight}, but any resolution
satisfying certain conditions will suffice to compute the bounded
cohomology in a completely canonical way.
More specifically, the prominent role played by proper actions
in the case of continuous cohomology is played by amenable actions
in the case of bounded continuous cohomology.

\begin{thm}[Burger--Monod \cite{Burger_Monod_GAFA}, Monod \cite{Monod_book}]\label{thm:amen}
Let $G$ be a locally compact second
countable group and $(S,\nu)$ a regular amenable $G$-space.
Then the continuous bounded cohomology of $G$ is isometrically
isomorphic to the cohomology of the complex
\bqn
\xymatrix@1{
 0\ar[r]
&\la(S,\RR)^G\ar[r]^-d
&\la(S^2,\RR)^G\ar[r]^-d
&\dots
}
\eqn
with the usual homogeneous coboundary operator.
\end{thm}
Here $\la(S^n,\RR)$ denotes the subspace of $\linfty(S^n,\RR)$
consisting of functions such that 
$f(s)=\operatorname{sign}(\sigma)f\big(\sigma(s)\big)$
for all $s\in S^n$ and $\sigma$ any permutation of the coordinates.

Without getting into the details of the amenability of an action
(for which we refer the reader to \cite{Zimmer_book}),
let us mention that the action of a group $\Lambda$
on the Poisson boundary $(B,\nu)$ relative to a probability measure $\theta$
is amenable, 
as well as the action of a connected semisimple Lie group $G$
on the quotient $G/P$ by a minimal parabolic subgroup $P<G$.
So, for example, the action of a surface group $\gg$ on $\SS^1$
via a hyperbolization is amenable, but if $\Xx$ is a Hermitian symmetric
space the action of $\Isom(\Xx)^\circ$ on the Shilov boundary $\cs_\Xx$ 
is not, unless the symmetric space has real rank one.

If in addition to being amenable the action of $G$ on $(S,\nu)$
is mixing, that is the diagonal action on $(S\times S,\nu\times \nu)$
is ergodic, then any $G$-invariant measurable function on $S\times S$ 
must be essentially constant, and hence $\la(S^2,\RR)^G=0$.
This, together with Theorem~\ref{thm:amen} implies the following

\begin{cor}\label{cor:mixing} Let $G$ be a locally compact second countable group
and $(S,\nu)$ a regular amenable mixing $G$-space. 
If $\Zz\la(S^3,\RR)$ denotes the subspace of cocycles in $\la(S^3,\RR)$,
then we have a canonical isometric isomorphism 
\bqn
\hcb^2(G,\RR)\cong\Zz\la(S^3,\RR)^G\,.
\eqn
\end{cor}

\begin{exo}\label{exo:coh} Since the $\gg$-action on $\SS^1$ is amenable and mixing,
then 
\bqn
\hb^2(\gg,\RR)\cong\Zz\la\big((\SS^1)^3,\RR\big)^{\gg}\,.
\eqn
Likewise if $G$ is a connected semisimple Lie group and $P<G$
is a minimal parabolic, then 
\bqn
\hcb^2(G,\RR)\cong\Zz\la\big((G/P)^3,\RR\big)^{G}\,.
\eqn
\end{exo}

On the one hand this shows immediately that in degree two
continuous bounded cohomology is a Banach space, on the other
it allows us to represent bounded cohomology classes via meaningful
cocycles defined on boundaries.  

\medskip

From now on we shall apply these considerations to the symplectic group
$G=\sp(V)$; for ease of notation, set $\dim V=2n$.
Following Kashiwara \cite[\S~1.5]{Lion_Vergne}, 
we recall that the {\it Maslov index}
$\beta_n$ of three Lagrangians $L_1,L_2,L_3\in\Ll(V)$ is defined as the index
$\beta_n(L_1,L_2,L_3)\in\ZZ$ of the quadratic form
\bqn
\ba
L_1\oplus L_2\oplus L_3&\longrightarrow\qquad\qquad\RR\\
(x_1,x_2,x_3)\,\,\,&\mapsto\<x_1,x_2\>+\<x_2,x_3\>+\<x_3,x_1\>\,.
\ea
\eqn
The function $\beta_n:\Ll(V)^3\to\ZZ$ is a cocycle 
which takes integer values in the interval $[-n,n]$; 
more specifically, on the space $\Ll(V)^{(3)}$ 
of triples of Lagrangians which are pairwise transverse, its set of values is
$\{-n,-n+2,\dots,n-2,n\}$, and each fiber of $\beta_n$ is precisely
an open $\sp(V)$-orbit.  Remark also that $\beta_1$ is nothing but 
the orientation cocycle on $\SS^1$.

The space $\Ff(V)$ of complete isotropic flags is a homogeneous space
of $\sp(V)$ with a minimal parabolic subgroup as stabilizer,
and therefore the $\sp(V)$-action on $\Ff(V)$ is amenable.  Let 
\bqn
\pr:\Ff(V)\to\Ll(V)
\eqn
be the projection
\bqn
\pr\big(\{0\}\subsetneq V_1\subsetneq\dots\subsetneq V_n\big):= V_n\,.
\eqn
With these notations we have:

\begin{prop}\label{prop:maslov} The map
\bq\label{eq:repr}
\beta_n\circ\pr^3:\Ff(V)^3\to\ZZ
\eq
is a bounded $\sp(V)$-invariant alternating cocycle such that
$\pi(\beta_n\circ\pr^3)$ corresponds to the bounded \kahler class
$\kappa_{\sp(V)}^\mathrm{b}\in\hcb^2(\sp(V,\RR)$ under
the isometric isomorphism in Corollary~\ref{cor:mixing}.
In particular
\bqn
\big\|\kappa_{\sp(V)}^\mathrm{b}\big\|=\|\pi(\beta_n\circ\pr^3)\|_\infty=\pi\,n\,.
\eqn
\end{prop}

Of course the drawback of the acquired freedom in going from continuous
functions to $\linfty$ functions -- or, more specifically, function classes -- 
is that now the implementation of the pullback of a bounded cohomology class 
cannot be done mindlessly as before,
since pullbacks even via continuous maps do not define,
in general, a well defined equivalence class of measurable functions.
However, the situation is much simpler in our case,
given that our class admits as a representative the Borel function in (\ref{eq:repr})
for which the cocycle identity holds everywhere.
The following important result is a particular case 
of a general phenomenon for which we refer the reader to \cite{Burger_Iozzi_app}.

\begin{thm}\label{thm:pullback}  Let $\gg\to\sp(V)$ be a homomorphism,
and assume that there exists a $\rho$-equivariant measurable map
\bqn
\varphi:\SS^1\to\Ll(V)\,,
\eqn
where $\gg$ acts on $\SS^1$ via a hyperbolization.
Then the pullback 
\bqn
\rho^{(2)}_\mathrm{b}\big(\kappa_{\sp(V)}^\mathrm{b}\big)\in\hcb^2(\gg,\RR)
\eqn
is represented by the cocycle $\pi(\beta_n\circ\varphi^3):(\SS^1)^3\to\RR$ defined by
\bqn
(x,y,z)\mapsto\pi\beta_n\big(\varphi(x),\varphi(y),\varphi(z)\big)\,.
\eqn
\end{thm}

Now we succeeded in implementing the pullback in a rather effective way,
but we find ourselves in the infinite dimensional Banach space
$\hb^2(\gg,\RR)$.  To size things down again, we shall need to make
use of the transfer map.

Choose a hyperbolization of $\sg$ and let as before $\Gamma$ be the 
realization of $\gg$ as a cocompact lattice in $\pu(1,1)$. Inspired
by Example~\ref{exo:coh} and by the fact that 
\bqn
\hcb^2\big(\pu(1,1),\RR\big)\cong\Zz\la\big((\SS^1)^3,\RR\big)^{\pu(1,1)}\,,
\eqn
define a transfer map
\bqn
t:\linfty\big((\SS^1)^3,\RR\big)^\Gamma\to\linfty\big((\SS^1)^3,\RR\big)^{\pu(1,1)}
\eqn
by
\bqn
tf(x,y,z):=\int_{\Gamma\backslash\pu(1,1)}f(gx,gy,gz)\,d\mu(g)\,,
\eqn
where $\mu$ is the $\pu(1,1)$-invariant probability measure on 
$\Gamma\backslash\pu(1,1)$.  Since by Proposition~\ref{prop:maslov}
\bqn
\hcb^2\big(\sp(V),\RR\big)\cong\RR\cdot(\beta_n\circ\pr^3)
\eqn
and
\bqn
\hcb^2\big(\pu(1,1),\RR\big)\cong\RR\cdot\beta_1\,,
\eqn
composition of the pullback implemented as in Theorem~\ref{thm:pullback}
followed by the transfer map in cohomology 
\bq\label{eq:compo}
\xymatrix@1{
 \hcb^2\big(\sp(V),\RR\big)\ar[r]^-{\rho^{(2)}_\mathrm{b}}
&\hb^2(\Gamma,\RR)\ar[r]^-{t^{(2)}}
&\hcb^2\big(\pu(1,1),\RR\big)
}
\eq
implies that there exists a constant $c\geq0$ such that for almost all
$x,y,z\in\SS^1$
\bq\label{eq:half-formula}
\int_{\Gamma\backslash\pu(1,1)}
   \beta_n\big(\varphi(gx),\varphi(gy),\varphi(gz)\big)\,d\mu(g)=
   c\beta_1(x,y,z)\,.
\eq
An analogous composition of maps as in (\ref{eq:compo})
in ordinary cohomology and their interplay via the comparison map 
which for $\sp(V)$ and $\pu(1,1)$ are isomorphisms \cite{Burger_Monod_GAFA}, 
allow us to explicit the constant $c$ in (\ref{eq:half-formula}) 
as explained in \cite[\S~3]{Iozzi_ern} in the context of Matsumoto's theorem.

\begin{thm}\label{thm:formula} Let $\rho:\gg\to\sp(V)$ be a homomorphism,
$\Gamma<\pu(1,1)$ a hyperbolization of $\gg$, and assume that there exists
a $\rho$-equivariant measurable map $\varphi:\SS^1\to\Ll(V)$.
Then for almost every $x,y,z\in \SS^1$
\begin{equation}
\int_{\Gamma\backslash\pu(1,1)}
   \beta_n\big(\varphi(gx),\varphi(gy),\varphi(gz)\big)\,d\mu(g)=
   \frac{\tr}{\big|\chi(\sg)\big|}\beta_1(x,y,z)\,.
\end{equation}
\end{thm}

Observe that if either $\rho(\Gamma)$ is Zariski dense or $\rho$
is tight, such a measurable $\Gamma$-equivariant map exists.
The following corollary is then immediate from Theorem~\ref{thm:map}
and Theorem~\ref{thm:formula}.

\begin{cor}\label{cor:formula} Let $\rho:\Gamma\to\sp(V)$ be a 
maximal representation.
Then there exists a $\rho$-equivariant measurable map
$\varphi:\SS^1\to\Ll(V)$ and it satisfies
\bqn
\beta_n\big(\varphi(x),\varphi(y),\varphi(z)\big)=n\beta_1(x,y,z)
\eqn
for almost every $x,y,z\in\SS^1$.
\end{cor}


\vskip1cm

\section{Symplectic Anosov Structures:  Proofs}\label{sec:anproofs}
In this section we prove the results stated in \S~\ref{sec:anosov}.
These proofs rest entirely on Corollary~\ref{cor:formula}
and are otherwise independent of the machinery used 
to establish Corollary~\ref{cor:formula}.

\subsection{The Geometry of Triples of Lagrangians}
Here we collect a few basic facts about the Maslov cocycle.
Our reference is \cite[\S~1.5]{Lion_Vergne}.

The space $\Ll(V)^{(3)}$ of triples of pairwise transverse
Lagrangians decomposes as a union $\sqcup_{j=0}^n\Oo_{n-2j}$
of $(n+1)$ open $\sp(V)$-orbits such that $\Oo_{n-2j}$ is 
the level set of $\beta_n$ where $\beta_n$ takes the value $n-2j$.

The maximal value $n$ is special in that, 
if $L_1, L_2,L_3$ are not pairwise transverse,
then $\big|\beta_n(L_1,L_2,L_3)\big|<n$, \cite[Proposition~1.5.10]{Lion_Vergne}.
Thus we observe that
\begin{equation}\label{propr:max-tr}
\begin{aligned}
&\hbox{if }\beta_n(L_1,L_2,L_3)=\pm n,\\
\hbox{then }\,L_1,&L_2,L_3\hbox{ are pairwise transverse.}
\end{aligned}
\end{equation}
Given $L_1,L$ and $L_3$ with $L_1$ and $L$ transverse to $L_3$,
consider the linear map $T_{13}:L_1\to L_3$ defined by 
\bqn
L=\big\{\ell_1+T_{13}(\ell_1):\,\ell_1\in L_1\big\}
\eqn
and the quadratic form $Q_{L}^{L_1,L_3}:L_1\to\RR$ defined by
\bqn
Q_{L}^{L_1,L_3}(x):=\big\<x, T_{13}x\big\>\,.
\eqn

\begin{figure}[!h]
\centering
\psfrag{l1}{$\ell_1$}
\psfrag{l3}{$\ell_3$}
\psfrag{l1+l3}{$\ell_1+\ell_3$}
\psfrag{L1}{$L_1$}
\psfrag{L2}{$L$}
\psfrag{L3}{$L_3$}
\psfrag{o}{$0$}
\includegraphics[width=1\linewidth]{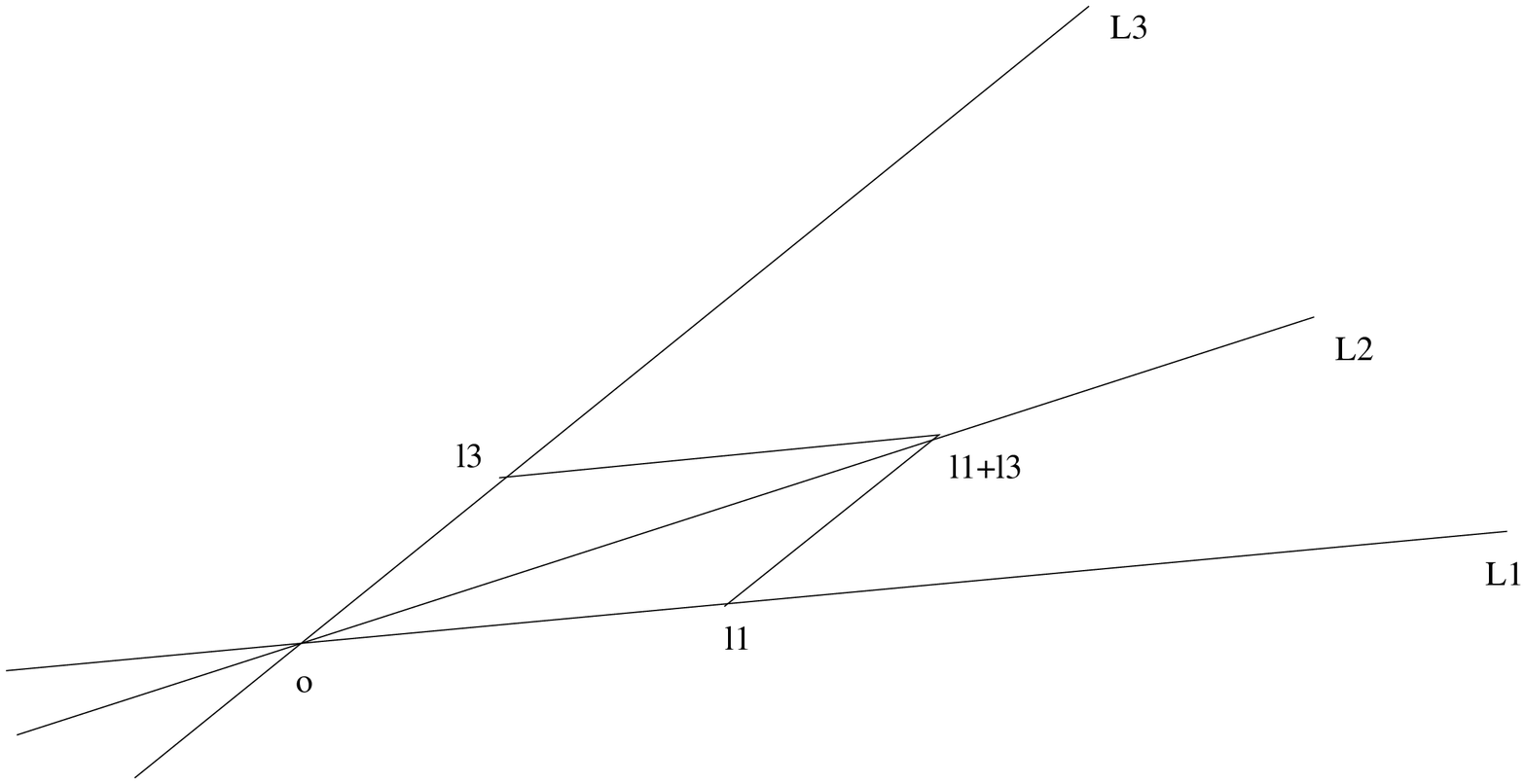}
\caption{\small Let $L_1$ and $L$ be transverse to $L_3$.  Then 
the vector $\ell_3\in L_3$ is the image of the vector $\ell_1\in L_1$
under the isomorphism $T_{13}:L_1\to L_3$ defined by $L$.
The value at $\ell_1$ of the quadratic form $Q_L^{L_1,L_3}$ on $L_1$ 
based at $L_3$ and induced by $L$ measures
the signed area of the parallelogram with vertices $0,\ell_1,\ell_1+\ell_3,\ell_3$.
Moreover, if also $L_1$ and $L$ are transverse, then 
$Q_L^{L_1,L_3}(\ell_1)=-Q_L^{L_3,L_1}(\ell_3)$, where $\ell_1$ and
$\ell_3$ are related as above.}
\label{fig:quadr}
\end{figure}

Let now
\bqn
t(L_3):=\big\{L\in\Ll(V):\,L\cap L_3=\{0\}\big\}
\eqn
and let $\Qq(L_1)$ be the space of quadratic forms on $L_1$.
Then we have a diffeomorphism
\bq\label{eq:diff}
\ba
t(L_3)&\to \Qq(L_1)\\
L\,\,\,&\mapsto Q_L^{L_1,L_3}
\ea
\eq
and moreover (see \cite[Lemma~1.5.4]{Lion_Vergne})
\bq\label{eq:sforms}
\beta_n(L_1,L,L_3)=\operatorname{sign}\big(Q_L^{L_1,L_3}\big)\,.
\eq
If $\tau:=(L_1,L_2,L_3)$ is a triple of pairwise transverse Lagrangians,
we have an endomorphism $J(\tau)$ of $V=L_1\oplus L_3$
given in block form by 
\bq\label{eq:complex-structure}
J(\tau):=\begin{pmatrix}0&-T_{31}\\T_{13}&0\end{pmatrix}
\eq
which, since $J(\tau)^2=-\id$, defines a complex structure on $V$;
moreover $\<\,\cdot\,,J(\tau)\cdot\,\>$ is symmetric and
the associated quadratic form $q_{J(\tau)}$ is the orthogonal direct sum
of $Q_{L_2}^{L_1,L_3}$ on $L_1$ and $-Q_{L_2}^{L_3,L_1}$ on $L_3$
(see Figure~\ref{fig:quadr});
in particular $q_{J(\tau)}$ has signature 
\bqn
2\beta_n(L_1,L_2,L_3)=\operatorname{sign}\big(Q_{L_2}^{L_1,L_3}\big)-
                     \operatorname{sign}\big(Q_{L_2}^{L_3,L_1}\big)\,.
\eqn

If now $\Ll(V)^3_{max}$ denotes the set of triples $\tau$
for which $\beta_n(\tau)=n$, we obtain 
an $\sp(V)$-equivariant map
\bq\label{eq:lvtox}
\ba
\Ll(V)^3_{max}\to&\,\,\Xx_{\sp(V)}\\
\tau\,\,\,\,\longmapsto& \,\,\,J(\tau)
\ea
\eq
into the symmetric space $\Xx_{\sp(V)}$ associated to $\sp(V)$.

\begin{defi} We say that a quadruple $\tau'$ of Lagrangians
is {\it maximal} if $\beta(\tau)=n$ for any subtriple
of Lagrangians $\tau$ taken in the same cyclic order as in $\tau'$.
\end{defi}

In particular (\ref{propr:max-tr}) implies that a maximal quadruple
consists of pairwise transverse Lagrangians.
Finally we have the following important monotonicity property:

\begin{lemma}\label{lem:monotone} 
Assume that the quadruple of Lagrangians $(L_0,L_1,L_2,L_\infty)$ is maximal. 
Then 
\bqn
0<Q_{L_1}^{L_0,L_\infty}<Q_{L_2}^{L_0,L_\infty}
\eqn 
and 
\bqn
Q_{L_1}^{L_\infty,L_0}<Q_{L_2}^{L_\infty,L_0}<0\,.
\eqn 
\end{lemma}
\begin{proof}  For $\ell_0\in L_0$, let $\ell_\infty,\ell_\infty'\in L_\infty$
with $\ell_0+\ell_\infty\in L_1$ and $\ell_0+\ell'_\infty\in L_2$. Then
\bqn
\ba
 Q_{L_2}^{L_0,L_\infty}(\ell_0)-Q_{L_1}^{L_0,L_\infty}(\ell_0)
&=\<\ell_0,\ell_\infty'-\ell_\infty\>\\
&=\<\ell_0+\ell_\infty,\ell_\infty'-\ell_\infty\>\\
&=Q_{L_2}^{L_1,L_\infty}(\ell_0+\ell_\infty)\,,
\ea
\eqn
where the last equality follows from the fact that $(\ell_0+\ell_\infty)\in L_1$,
$\ell_\infty'-\ell_\infty\in L_\infty$, and their sum $\ell_0+\ell_\infty'\in L_2$.
Maximality of $(L_1, L_2,L_\infty)$ implies that 
$Q_{L_2}^{L_1,L_\infty}>0$, and 
maximality of $(L_0,L_1,L_\infty)$ implies that $Q_{L_1}^{L_0,L_\infty}>0$.
Hence the assertion.
\end{proof}

Notice that in the proof of Lemma~\ref{lem:monotone} what was used is
exactly the fact that 
the Lagrangians $L_0,L_1,L_2,L_\infty$ are pairwise transverse 
and that the triples $(L_0,L_1,L_\infty)$ and $(L_1,L_2,L_\infty)$ are maximal,
which however, via the cocycle identity for $\beta_n$, 
is equivalent to the maximality of the quadruple $(L_0,L_1,L_2,L_\infty)$
(see the proof of Lemma~\ref{lem:transverse}).

\subsection{Proofs of the Results in \S~\ref{sec:anosov}}
Let $\rho:\Gamma\to\sp(V)$ be a maximal representation 
and let $\varphi:\SS^1\to\Ll(V)$ be the $\rho$-equivariant
measurable map given by Corollary~\ref{cor:formula}.
Paramount in the study of regularity properties of the map $\varphi$
is the closer analysis of its essential graph which we now introduce.
Let $\lambda$ be the Lebesgue measure on $\SS^1$.  The essential graph
$\Ef$ of $\varphi$ is the closed subset $\Ef\subset\SS^1\times\Ll(V)$ 
which is the support of the pushforward of the measure $\lambda$ under
the map
\bqn
\ba
\SS^1&\to\SS^1\times\Ll(V)\\
x&\mapsto\big(x,\varphi(x)\big)\,.
\ea
\eqn
Here and in the sequel we shall often use the observation 
that 
\bq\label{eq:essgr}
\hbox{ for almost every }x\in\SS^1,\,\,\big(x,\varphi(x)\big)\in\Ef\,.
\eq

\begin{lemma}\label{lem:valessgr}
Let $(x_1,L_1),(x_2,L_2),(x_3,L_3)\in\Ef$, and assume that
\be
\item $x_1,x_2,x_3$ are pairwise distinct, and
\item $L_1,L_2,L_3$ are pairwise transverse.
\ee
Then $\beta_n(L_1,L_2,L_3)=n\beta_1(x_1,x_2,x_3)$\,.
\end{lemma}
\begin{proof} We may assume that $\beta_1(x_1,x_2,x_3)=1$.
Using that $(x_i,L_i)\in\Ef$, Corollary~\ref{cor:formula} 
and the definition of essential graph imply that
we may find sequences $L_i^{(k)}$, $i=1,2,3$, $k\in\NN$, such that 
$\beta_n\big(L_1^{(k)},L_2^{(k)},L_3^{(k)}\big)=n$ and 
$\big(L_1^{(k)},L_2^{(k)},L_3^{(k)}\big)$ converges to $(L_1,L_2,L_3)$.
In particular, $(L_1,L_2,L_3)$ is in the closure 
$\overline{\Oo_n}$ in $\Ll(V)^3$ of $\Oo_n$.  Since 
on the other hand this triple belongs to $\sqcup_{j=0}^n\Oo_{n-2j}$,
observing that $\Oo_k\cap\overline{\Oo_n}=\emptyset$ for $k\neq n$,
we conclude that $(L_1,L_2,L_3)\in\Oo_n$.  
\end{proof}

Notice now that any two (distinct) points $x_1,x_2\in\SS^1$ 
determine an interval in $\SS^1$, by defining 
\bqn
((x_1,x_2)):=\{t\in\SS^1:\,\beta_1(x_1,t,x_2)=1\}\,.
\eqn

\begin{lemma}\label{lem:transverse} Let $(x_1,L_1)$ and $(x_2,L_2)\in\Ef$ with $x_1\neq x_2$.
Then $L_1$ and $L_2$ are transverse.
\end{lemma}
\begin{proof}
Using Corollary~\ref{cor:formula}, (\ref{eq:essgr}) and Remark~\ref{rem:way-out}
twice, we may choose
$a\in((x_1,x_2))$ such that $\big(a,\varphi(a)\big)\in\Ef$ and
$\varphi(a)$ is transverse to $L_1,L_2$, and choose $b\in((x_2,x_1))$
so that $\big(b,\varphi(b)\big)\in\Ef$ and
$\varphi(b)$ is transverse to $\varphi(a),L_1,L_2$.

Applying the cocycle property of $\beta_n$ to the quadruple
$\varphi(a),L_2,\varphi(b),L_1$, we have that
\bqn
\ba
&\beta_n\big(L_2,\varphi(b),L_1\big)-\beta_n\big(\varphi(a),\varphi(b),L_1\big)\\
+&\beta_n\big(\varphi(a),L_2,L_1\big)-\beta_n\big(\varphi(a),L_2,\varphi(b)\big)=0\,;
\ea
\eqn
since it follows from Lemma~\ref{lem:valessgr} that 
\bqn
\beta_n\big(\varphi(a),\varphi(b),L_1\big)=n=
\beta_n\big(\varphi(a),L_2,\varphi(b)\big)\,,
\eqn
we obtain that 
\bqn
\beta_n\big(L_2,\varphi(b),L_1\big)+\beta_n\big(\varphi(a),L_2,L_1\big)=2n\,,
\eqn
which implies in turn that
\bqn
\beta_n\big(L_2,\varphi(b),L_1\big)=
\beta_n\big(\varphi(a),L_2,L_1\big)=n\,.
\eqn
It follows hence from (\ref{propr:max-tr}) that $L_1$ and $L_2$ are transverse.
\end{proof}

From Lemmas~\ref{lem:valessgr} and \ref{lem:transverse} we deduce the following

\begin{cor}\label{cor:max} For $(x_1,L_1),(x_2,L_2),(x_3,L_3)\in\Ef$
with $(x_1,x_2,x_3)$ pairwise distinct, we have
\bqn
\beta_n(L_1,L_2,L_3)=n\beta_1(x_1,x_2,x_3)\,.
\eqn
\end{cor}

For the following, it will be convenient to define for $A\subset\SS^1$
the ``image of $A$'' by $\Ef$
\bqn
F_A:=\big\{L\in\Ll(V):\,\hbox{ there exists }a\in A\hbox{ such that }(a,L)\in\Ef\big\}
\eqn
which is closed if $A\subset\SS^1$ is so. Now let us fix  any two distinct
points $x,y\in\SS^1$.

\begin{lemma}\label{lem:sets}
The sets $\overline F_{((y,x))}\cap F_{\{x\}}$ and
$\overline F_{((x,y))}\cap F_{\{x\}}$ both consist of one point.
\end{lemma}
\begin{proof}  Assume that there are $L_0,L'_0\in\overline F_{((x,y))}\cap F_{\{x\}}$
and fix $L_\infty\in F_{\{y\}}$.  By hypothesis, there are
sequences $(x_n,L_n)$ and $(x'_n,L'_n)$ in $\Ef$ with
\be
\item $x_n,x'_n\in((x,y))$, and $\lim x_n=\lim x'_n=x$;
\item $\lim L_n=L_0$ and $\lim L'_n=L'_0$\,.
\ee
By Lemma~\ref{lem:transverse} all $L_n$ and $L'_n$ are transverse to $L_\infty$
and we may thus use the diffeomorphism in (\ref{eq:diff})
\bqn
\ba
t(L_\infty)&\to \Qq(L_0)\\
L\quad&\mapsto Q_L^{L_0,L_\infty}
\ea
\eqn
and study the situation in the model $\Qq(L_0)$.
Dropping the superscript $L_0,L_\infty$, we have 
that $\lim Q_{L_n}=Q_{L_0}=0$.  For every $k\geq1$,
there is $N(k)$ such that $x'_n\in((x,x_k))$ for all $n\geq N(k)$,
and consequently $L_0,L'_n,L_k,L_\infty$ is maximal;
using Lemma~\ref{lem:monotone}, this implies that 
\bqn
Q_{L_0}=0\leq Q_{L'_n}\leq Q_{L_k}
\eqn
and hence $\lim_nQ_{L'_n}=0$.  This shows that 
$\lim_n L'_n=L_0$ and hence $L'_0=L_0$.
\end{proof}

According to Lemma~\ref{lem:sets}, for every $x\in\SS^1$
define
\bqn
\ba
\varphi_+(x)\in&\overline F_{((y,x))}\cap F_{\{x\}}\hbox{ and }
\varphi_-(x)\in&\overline F_{((x,y))}\cap F_{\{x\}}\,.
\ea
\eqn
From the definitions one deduces immediately the following

\begin{cor}\label{cor:fi-fi+}
The maps 
\bqn
\varphi_+,\varphi_-:\SS^1\to\Ll(V)
\eqn
defined above are respectively left and right continuous
and strictly $\Gamma$-equivariant.
\end{cor}

\bigskip

\begin{figure}[!h]
\centering
\psfrag{u-}{$u_-$}
\psfrag{u0}{$u_0$}
\psfrag{ut}{$u_t$}
\psfrag{u+}{$u_+$}
\psfrag{u}{$u$}
\psfrag{gtu}{$g_tu$}
\includegraphics[width=.5\linewidth]{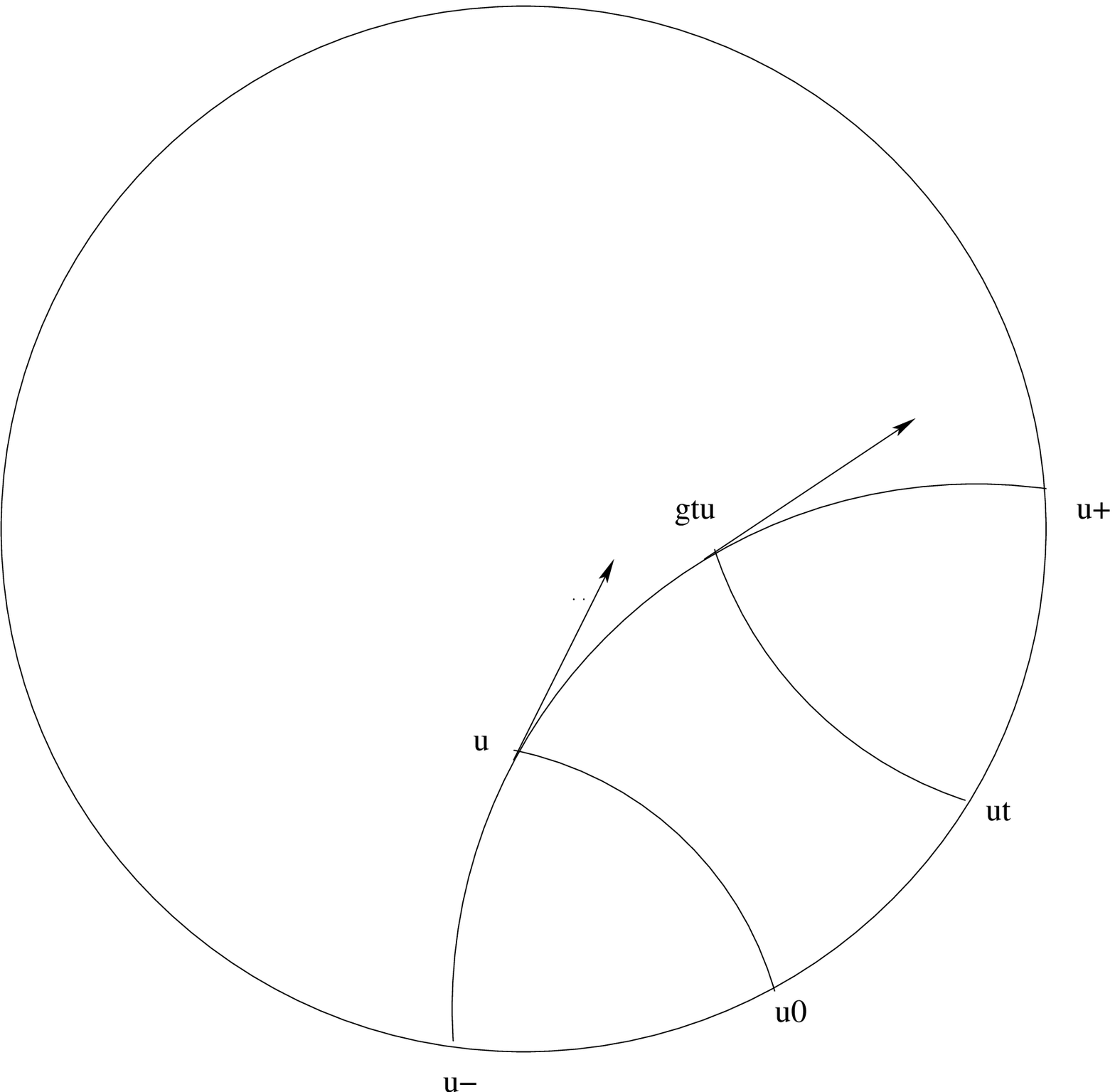}
\caption{\small {The identification of $T^1\Dd_{1,1}$ with $(\SS^1)^{(3)}$.}}
\label{fig:tgbdl}
\end{figure}

Now we turn to our symplectic bundle $E^\rho$ introduced 
in \S~\ref{sec:anosov} and the study of the properties of the flow $g_t^\rho$.
To define the Lagrangian splitting of $E^\rho$ we parametrize
$T^1\Dd_{1,1}$ by the set $(\SS^1)^{(3)}$ of distinct triples of points on $\SS^1$,
as follows: to a unit vector $u\in T^1\Dd_{1,1}$
based at $x$ associate the triple $(u_-,u_0,u_+)\in\SS^1$, where
$u_-\in\SS^1$ and $u_+\in\SS^1$ are respectively the initial and ending point
of the geodesic $[u_-,u_+]$ determined by $u$,
and $u_0\in\SS^1$ is the endpoint of the geodesic perpendicular to 
$[u_-,u_+]$ at $x\in\Dd_{1,1}$ and oriented in such a way that $u_0\in((u_-,u_+))$.
Notice that as $u$ moves along the geodesic $[u_-,u_+]$ in the positive
direction, the point $u_0$ approaches $u_+$ but the points $u_-,u_+$
stay unchanged, so that the vector $g_tu$ corresponds to the triple
$(u_-,u_t,u_+)$ (see Figure~\ref{fig:tgbdl}).

Let $\varphi_-,\varphi_+:\SS^1\to\Ll(V)$ be respectively the right
and left continuous $\Gamma$-equivariant map in Corollary~\ref{cor:fi-fi+}.
For every $u\in T^1\Dd_{1,1}$, since $u_-\neq u_+$, 
Lemma~\ref{lem:transverse} implies that $\varphi_-(u_-)$ and $\varphi_+(u_+)$ 
define transverse and hence complementary Lagrangians
\bqn
V=\varphi_-(u_-)\oplus\varphi_+(u_+)\,.
\eqn
In this way we obtain a splitting of $\widetilde E^\rho$
into $\big(\tilde g_t^\rho\big)$-invariant Borel subbundles
$\widetilde E^\rho=\widetilde E_-^\rho\oplus\widetilde E_+^\rho$
which descends to a $\big(g_t^\rho\big)$-invariant splitting
\bqn
E^\rho=E_-^\rho\oplus E_+^\rho\,.
\eqn  
Using Corollary~\ref{cor:max} we deduce that the triples 
$\big(\varphi_-(u_-),\varphi_\pm(u_t),\varphi_+(u_+)\big)$
are maximal for every $t$, so that we can associate to each of them
complex structures $J(g_tu,+)$ and $J(g_tu,-)$ on $V$
as in (\ref{eq:lvtox}), and hence positive quadratic forms
$q_{J(g_tu,+)}$ and $q_{J(g_tu,-)}$, which thus give rise
to two families $\|\cdot\|_{g_tu}^+$ and $\|\cdot\|_{g_tu}^-$
of Euclidean metrics on $E^\rho(g_tu)$, for $t\in\RR$, $u\in T^1\Dd_{1,1}$.

\begin{lemma}\label{lem:mess}
Let $p:E^\rho\to T^1\Sigma$ be the projection defined in (\ref{eq:p}) 
and, if $\xi\in E^\rho$, let $u:=p(\xi)\in T^1\Sigma$.  Then
\be
\item For every $\xi\in E_+^\rho$ 
\bqn
\lim_{t\to+\infty}\|g_t^\rho\xi\|_{g_tu}^+=0\hbox{ monotonically, and }
\|g_{-t}^\rho\xi\|^\pm_{g_{-t}u}\geq\|\xi\|^\pm_u\hbox{ for all }t\geq0\,.\eqn
\item For every $\xi\in E_-^\rho$ 
\bqn
\lim_{t\to+\infty}\|g_{-t}^\rho\xi\|_{g_{-t}u}^-=0\hbox{ monotonically, and }
\|g_t^\rho\xi\|^\pm_{g_tu}\geq\|\xi\|^\pm_u,\hbox{ for all }t\geq0\,.
\eqn
\ee
\end{lemma}

\begin{figure}[!h]
\centering
\psfrag{f-u-}{$\varphi_-(u_-)$}
\psfrag{f+u0}{$\varphi_+(u_0)$}
\psfrag{f+ut}{$\varphi_+(u_t)$}
\psfrag{f+u+}{$\varphi_+(u_+)$}
\includegraphics[width=1\linewidth]{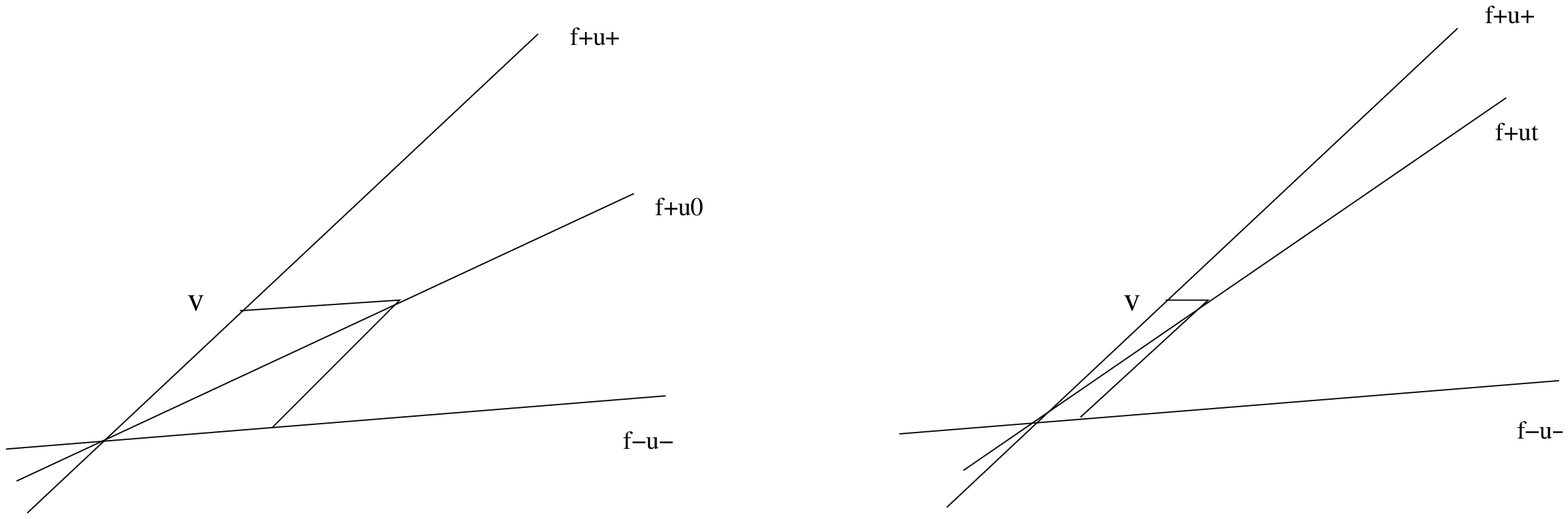}
\caption{\small {}}
\label{fig:limit}
\end{figure}

\begin{proof} We prove (1), as the proof of (2) is analogous.
Working in $\widetilde E^\rho$ as we may, 
let $\xi\in \widetilde E^\rho$, $\xi=(u,v)$, $v\in V$.  
Let $v\in\varphi_+(u_+)$.  We use the Euclidean metrics
$\|\cdot\|^+_{g_tu}$ defined by the triple
\bqn
\big(\varphi_-(u_-),\varphi_+(u_t),\varphi_+(u_+)\big)\,,
\eqn
that is
\bqn
 \big\|\tilde g_t^\rho\xi\big\|_{g_tu}^+
=\left|Q_{\varphi_+(u_t)}^{\varphi_+(u_+),\varphi_-(u_-)}(v)\right|\,,
\eqn
which, since $\varphi_+$ is left continuous and hence
\bqn
\lim_{t\to+\infty}\varphi_+(u_t)=\varphi_+(u_+)\,,
\eqn
implies immediately that 
\bqn
\lim_{t\to+\infty}\big\|\tilde g_t^\rho\xi\big\|^+_{g_tu}=0\,.
\eqn
Monotonicity follows from Lemma~\ref{lem:monotone}.  In fact,
for every $0\leq t_1<t_2$, the quadruple
\bqn
\big(\varphi_-(u_-),\varphi_+(u_{t_1}),\varphi_+(u_{t_2}),\varphi_+(u_+)\big)
\eqn
is maximal and hence Lemma~\ref{lem:monotone} implies that
\bqn
\big\|\tilde g_{t_2}^\rho\xi\big\|^+_{g_{t_2}u}\leq\big\|\tilde g_{t_1}^\rho\xi\big\|_{g_{t_1}u}^+\,.
\eqn
To prove the second statement in (1), observe that for $t\geq0$ the quadruple 
\bqn
\big(\varphi_-(u_-),\varphi_+(u_{-t}),\varphi_+(u_0),\varphi_+(u_+)\big)
\eqn
is maximal and hence Lemma~\ref{lem:monotone} implies that
\bqn
    \big\|\tilde g_{-t}^\rho\xi\big\|^+_{g_{-t}u}
   =\big|Q_{\varphi_+(u_{-t})}^{\varphi_+(u_+),\varphi_-(u_-)}(v)\big|
\geq
    \big|Q_{\varphi_+(u_0)}^{\varphi_+(u_+),\varphi_-(u_-)}(v)\big|
   =\|\xi\|_u^+\,.
\eqn
The statement for the metrics $\|\cdot\|^-_{g_tu}$ follows analogously.
\end{proof}

The metrics $\|\cdot\|_u^+$ and $\|\cdot\|_u^-$ are Borel metrics 
on the bundle $E^\rho$.  Since the basis
$T^1\Sigma$ is compact, any two continuous Euclidean metrics
on $E^\rho$ are equivalent: we have then

\begin{lemma}\label{lem:eq-metric} The metrics $\|\cdot\|_u^+$ and $\|\cdot\|_u^-$ 
are equivalent to a continuous metric. 
\end{lemma}

This follows easily from the following two facts:
\be
\item[-] The proper action of $\Gamma$ on $(\SS^1)^{(3)}$ has compact quotient.
\item[-] For any compact subset $C\subset(\SS^1)^{(3)}$, 
the set of metrics 
\bqn
\big\{\|\cdot\|_u^\pm:\,(u_-,u_0,u_+)\in C\big\}
\eqn
is bounded.
\ee

\begin{proof}[Proof of Theorem~\ref{thm:anosov}]
Fix a continuous Euclidean metric $\|\cdot\|$ on $E^\rho$.
Then it follows from Lemmas~\ref{lem:eq-metric} and \ref{lem:mess}
that 
\bqn
E_\pm^\rho:=\big\{\xi\in E^\rho:\lim_{t\to\pm\infty}\|g_t^\rho\xi\|=0\big\}\,.
\eqn
This implies by the following classical argument that 
the subbundles $E_+^\rho$ and $E_-^\rho$ are continuous.
Let $u_m$ be a converging sequence in $T^1\Sigma$ with limit $u$,
and let $F\subset E^\rho(u)$ be any accumulation point of
\bqn
\big\{E_+^\rho(u_m):m\geq1\big\}
\eqn 
in the Grassmann $n$-bundle of $E^\rho$.  
Let $\{m_k\}$ be a subsequence with $\lim_{k\to\infty}E_+^\rho(u_{m_k})=F$.
For every $\xi\in F$ take $\xi_k\in E_+^\rho(u_{m_k})$
with $\lim_{k\to\infty}\xi_k=\xi$.  Then the function
\bqn
\ba
\RR^+&\to\,\,\,\,\RR^+\\
t\,\,\,&\mapsto\big\|g_t^\rho\xi\big\|
\ea
\eqn
being a uniform limit on compacts of the sequence of functions
$t\mapsto\big\|g_t^\rho\xi_k\big\|$ which vanish at infinity,
vanishes at infinity as well, which implies that $\xi\in E_+^\rho(u)$
and hence $F\subseteq E_+^\rho(u)$; since both spaces have the same
dimension, we conclude that $F=E_+^\rho(u)$.  
This shows continuity of the splitting.

This implies by the definition of $\widetilde E_\pm^\rho$
that both maps $\varphi_+$ and $\varphi_-$ from $\SS^1$ to $\Ll(V)$
are continuous.  But this implies easily that $\varphi_-=\varphi_+$;
we shall denote from now on by $\varphi$ this continuous
$\Gamma$-equivariant map.  This implies now the first assertion of 
Corollary~\ref{cor:rect}.

We are thus in the following situation:  for every $u\in T^1\Dd_{1,1}$,
we have the splitting
\bqn
V=\varphi(u_-)\oplus\varphi(u_+),\quad u=(u_-,u_0,u_+)
\eqn
which gives rise to the splittings
\bqn
\widetilde E^\rho=\widetilde E^\rho_-\oplus\widetilde E^\rho_+
\eqn
and
\bqn
E^\rho=E^\rho_-\oplus E^\rho_+
\eqn
into continuous $\tilde g_t^\rho$ and $g_t^\rho$ invariant subbundles.
We denote by $J(u)\in \Xx_{\sp(V)}$ the complex structure associated to the triple
\bqn
\big(\varphi(u_-),\varphi(u_0),\varphi(u_+)\big)
\eqn
as in (\ref{eq:lvtox}).  It is now immediate that the map
\bq\label{eq:map}
\ba
T^1\Dd_{1,1}\to&\,\Xx_{\sp(V)}\\
u\quad\longmapsto&\, J(u)
\ea
\eq
gives a positive complex structure $J$ of $E^\rho$ with the required properties
(see (\ref{eq:complex-structure})).
Let $\|\cdot\|_u$ be the Euclidean metric on $E^\rho$ induced
by the quadratic form $q_{J(u)}$.

In the notation of Lemma~\ref{lem:mess}, 
we have $\|\cdot\|_u^+=\|\cdot\|_u^-=\|\cdot\|_u$
and hence for every $\xi\in E_\pm^\rho$ with $p(\xi)=u$
\bqn
\lim_{t\to\infty}\big\|g_{\pm t}\xi\big\|_{g_{\pm t}u}=0\quad\hbox{monotonically.}
\eqn
We claim now that there exists $T>0$ such that for every
$\xi\in E_+^\rho$,
\bqn
\big\|g_t^\rho\xi\big\|_{g_tu}\leq\frac12\|\xi\|_u\quad\hbox{ for }t\geq T\,.
\eqn
Indeed, if this were not the case, by Lemma~\ref{lem:mess} 
there would exist a sequence $\xi_n\in E_+^\rho$ 
and $T_n\to+\infty$ with $\|\xi_n\|=1$ 
and $\|g_{T_n}^\rho\xi_n\|_{g_{T_n}u_n}=\frac12$. 
We may assume that $\xi_n$ converges to a point $\xi\in E_+^\rho$.
Then the sequence of functions
\bqn
\ba
\RR^+&\longrightarrow\,\,\RR^+\\
t\,\,\,&\mapsto\|g_t^\rho\xi_n\|_{g_tu_n}
\ea
\eqn
converges uniformly on compact sets to
\bqn
t\mapsto\|g_t^\rho\xi\|_{g_tu}\,.
\eqn
But, by monotonicity, we have that
\bqn
\|g_t^\rho\xi_n\|_{g_tu_n}\geq\frac12,\quad\hbox{ for } t\in[0,T_n],
\eqn
and since $T_n\to+\infty$, we deduce that 
\bqn
\|g_t^\rho\xi\|_{g_tu}\geq\frac12\quad\hbox{ for all }t\geq0\,,
\eqn
which contradicts the fact that $\xi\in E_+^\rho$.
Applying the inequality
\bqn
\big\|g_T^\rho\xi\big\|\leq\frac12\|\xi\|
\eqn
to $nT$, for $n\in \NN$, we obtain the exponential decay.
\end{proof}

\begin{proof}[Proof of Corollary~\ref{cor:qi}]
The proof will rely on the metric properties of the map
defined in (\ref{eq:map}).

Fix a unit tangent vector $v\in T^1\Dd_{1,1}$ based at $0\in\Dd_{1,1}$ and
let $J_0:=J(v)\in\Xx_{\sp(V)}$ be the corresponding complex structure on $V$.
Observe first of all that 
$d\big(J_0,\rho(\gamma)J_0\big)$ is bounded above linearly by the word length
$\ell(\gamma)$ of $\gamma$, as an argument by recurrence on $\ell(\gamma)$
easily shows.  In order to show the lower bound, we shall use the
contraction--dilation property of the Anosov flow in Theorem~\ref{thm:anosov}(2).

The essential step is estimating the distance in $\Xx_{\sp(V)}$
between $J(u)$ and $J(g_tu)$, 
given by 
\bqn
d\big(J(u),J(g_tu)\big)
=\big|\ln\|\id\|_{J(u),J(g_tu)}\big|+\big|\ln\|\id\|_{J(g_tu),J(u)}\big|\,,
\eqn
for any $u\in T^1\Dd_{1,1}$ and any $t\geq0$ (see \S~\ref{sec:anosov}). 

For $x\in\varphi(u_-)$, applying Theorem~\ref{thm:anosov},
we have that 
\bqn
q_{J(g_tu)}(x)\geq e^{2At}q_{J(u)}(x)
\eqn
and likewise for $x\in \varphi(u_+)$
\bqn
q_{J(g_tu)}(x)\leq e^{-2At}q_{J(u)}(x)\,.
\eqn
These inequalities, together with the fact that $\varphi(u_-)\oplus\varphi(u_+)$
is an orthogonal decomposition for both $q_{J(u)}$ and $q_{J(g_tu)}$,
imply that
\bqn
\|\id\|_{J(u),J(g_tu)}\geq e^{At}
\eqn
and
\bqn
\|\id\|_{J(g_tu),J(u)}\geq e^{At}\,,
\eqn
from which we deduce that
\bq\label{eq:lb}
d\big(J(u),J(g_tu)\big)\geq 2At\,.
\eq

Let now $\gamma\in\Gamma$ and let us choose $u\in T^1\Dd_{1,1}$ 
to be the tangent vector at $0\in\Dd_{1,1}$ to the geodesic segment
connecting $0$ to $\gamma0$ and let $t=d(0,\gamma0)$.
Applying (\ref{eq:lb}) to this situation and observing that
$g_tu=\gamma u$, we get that 
\bqn
d\big(J(u),\rho(\gamma)J(u)\big)\geq2 A d(0,\gamma 0)
\eqn
and hence 
\bqn
d\big(J_0,\rho(\gamma)J_0\big)\geq 2Ad(0,\gamma 0)-2C\,,
\eqn
where
\bqn
C:=\sup\big\{d\big(J(w_1),J(w_2)\big):\,w_1,w_2\hbox{ are based at }0\big\}\,.
\eqn
Finally, $d(0,\gamma 0)$ is bounded linearly below in terms
of $\ell(\gamma)$, as follows from the Milnor--Svarc lemma.
\end{proof}
%

\medskip
\noindent
{\it Proof of Corollary~\ref{cor:rect}.}\,\,  
The injectivity of the $\Gamma$-equivariant continuous map 
\bqn
\varphi:\SS^1\to\Ll(V)
\eqn
obtained in the proof of Theorem~\ref{thm:anosov},
follows for instance from Corollary~\ref{cor:formula} because of continuity.
So we finally turn to the proof of the rectifiability of the image of $\varphi$.
For this we shall put to use the $\sp(V)$-invariant causal structure on $\Ll(V)$.  

Let us fix $a\neq b\in\SS^1$, let $L_0:=\varphi(a)$ and $L_\infty:=\varphi(b)$,
so that on $\SS^1\setminus\{b\}$, $\varphi$ takes values in 
$t(L_\infty)$.  Composing the restriction of $\varphi$ to 
$\SS^1\setminus\{b\}$ with the usual diffeomorphism 
\bqn
\ba
t(L_\infty)&\to\Qq(L_0)\\
L\quad&\mapsto Q_L^{L_0,L_\infty}\,,
\ea
\eqn
gives rise to a continuous map
\bqn
c:\SS^1\setminus\{b\}\to\Qq(L_0)
\eqn
whose restriction to the interval $((a,b))$ has the 
following properties:
\be
\item it takes values in the cone $\Qq^+(L_0)$ of positive definite quadratic forms,
and
\item for every $t_1,t_2\in((a,b))$ such that $a,t_1,t_2,b$ are 
in positive cyclic order, $c(t_2)-c(t_1)\in\Qq^+(L_0)$.
\ee
Fixing a scalar product on $L_0$, we can identify $\Qq(L_0)$ with the space
$\operatorname{Sym}(L_0)$ of symmetric endomorphisms of $L_0$
and $\Qq^+(L_0)$ with the cone $\operatorname{Sym}^+(L_0)$ of positive definite
ones.  On $\operatorname{Sym}(L_0)$ we have a natural scalar product
\bqn
\<\<A,B\>\>:=\operatorname{tr} AB
\eqn
and we have that for every $A,B\in\operatorname{Sym}^+(L_0)$
\bqn
\<\<A,B\>\>>0\,,
\eqn
that is $\operatorname{Sym}^+(L_0)$ is an open convex acute cone.  
The assertion then follows from the following general fact

\begin{lemma} Let $C\subset E$ be an open convex acute cone in 
an Euclidean space and let $f:[0,1]\to C$ be a continuous
map such that for every $t_1<t_2$,
\bqn
f(t_2)-f(t_1)\in C\,.
\eqn
Then $f$ is of finite length.
\end{lemma}

\begin{proof} Fix $e\in C$.  We claim that since $C$ is acute
\bqn
k:=\inf_{x\in C}\frac{\<\<x,e\>\>}{\|x\|}>0\,.
\eqn
Indeed, otherwise there is a nonzero $x\in\overline C$ 
such that $\<\<x,e\>\>=0$.  
On the other hand, since $C$ is open, for $s<0$
and $|s|$ small enough we have that 
\bqn
e':=sx+(1-s)e\in C\,,
\eqn
which implies that $\<\<e',x\>\><0$
and contradicts the fact that $\<\<u,v\>\>\geq0$
for all $u,v\in\overline C$.  

Let $0\leq s<t\leq1$;  then $f(t)-f(s)\in C$
and applying the claim, we obtain:
\bqn
\|f(t)-f(s)\|\leq\frac{1}{k}\<f(t)-f(s),e\>\,.
\eqn
Given any subdivision $0=t_0<t_1<\dots<t_{n-1}<t_n=1$
of the interval $[0,1]$, we deduce that
\bqn
    \sum_{i=1}^n\|f(t_i)-f(t_{i-1})\|
\leq\frac{1}{k}\sum_{i=1}^n\<f(t_i)-f(t_{i-1}),e\>
   =\frac{\<f(1)-f(0),e\>}{k}
\eqn
which proves that $f$ is rectifiable.
\end{proof}


\vskip1cm

\bibliographystyle{amsplain}
\bibliography{refs}
\vskip1cm

\providecommand{\bysame}{\leavevmode\hbox to3em{\hrulefill}\thinspace}
\providecommand{\MR}{\relax\ifhmode\unskip\space\fi MR }
\providecommand{\MRhref}[2]{%
  \href{http://www.ams.org/mathscinet-getitem?mr=#1}{#2}
}
\providecommand{\href}[2]{#2}

\vskip1cm
\end{document}